\documentclass[11pt]{article}

\usepackage[utf8]{inputenc}
\usepackage{amsmath,amsfonts,latexsym, xspace,amsthm,graphicx, amssymb}
\usepackage{amsthm,mathrsfs,fancyhdr,subfigure}
\usepackage{graphicx} 
\usepackage[noadjust,sort]{cite}
\usepackage{enumerate}
\usepackage[letterpaper,left=1in,right=1in,top=1in,bottom=1in]{geometry}
\usepackage{hyperref}
\usepackage{authblk}

\usepackage{color,comment}
\usepackage{needspace}

\makeatletter
\renewcommand*{\@biblabel}[1]{\hfill#1.}
\makeatother

\makeatletter
\renewcommand{\fnum@figure}{Fig.\ \thefigure}
\makeatother

\newcommand{\Gnp}{G_{n,p}}
\newcommand{\Gncntext}{G_{n,c/n}}
\newcommand{\Gnm}{G_{n,m}}

\newcommand{\E}{{\mathbb E}}
\newcommand{\R}{{\mathbb R}}

\newcommand{\vol}{{\rm vol}}

\newcommand{\eps}{\varepsilon}

\newcommand{\cB}{\mathcal{B}}
\newcommand{\cC}{\mathcal{C}}
\newcommand{\matr}[1]{\mathbf{#1}}

\newcommand{\q}{q^*}
\newcommand{\qk}{q^*_{\leq k}}

\numberwithin{equation}{section}

\newtheorem{thm}{Theorem}[section]

\usepackage{rotating}

\let\stdcaption\caption
\usepackage{floatpag}
\let\caption\stdcaption

\usepackage{tikz}
\usepackage{tkz-berge}
\usepackage{tkz-graph}
\usepackage{tkz-base}
\usepackage{tkz-euclide}
\usepackage{pgf}

\newcommand{\pr}{{\mathbb P}}

\newcommand{\N}{{\mathbb N}}

\newcommand{\cA}{{\mathcal A}}
\newcommand{\cE}{{\mathcal E}}

\renewcommand{\eps}{\varepsilon}

\begin{document}

\title{Modularity and random graphs}

\author[1]{Colin McDiarmid\thanks{Email: \textit{cmcd@stats.ox.ac.uk}.}}
\author[2]{Fiona Skerman\thanks{Email: \textit{fiona.skerman@math.uu.se}. 
}}
\affil[1]{Department of Statistics, University of Oxford.}
\affil[2]{Department of Mathematics, Uppsala University.}

\date{\vspace{-20mm}}

\maketitle


\vspace*{.2in}

\begin{quote}\emph{$\!\!$
This work will appear as a chapter in a forthcoming volume titled `Topics in Probabilistic Graph Theory'. 
For a given graph $G$, each partition of the vertices has a modularity score, with higher values indicating that the partition better captures community structure in $G$. The modularity $\q(G)$ of $G$ is the maximum over all vertex-partitions of the modularity score, and satisfies $0\leq \q(G)< 1$. Modularity lies at the heart of the most popular algorithms for community detection.
In this chapter we discuss the behaviour of the modularity of various kinds of random graphs, starting with the binomial random graph~$\Gnp$ with $n$ vertices and edge-probability $p$.}
\end{quote}

\vspace*{.3in}

\noindent 1. \hspace*{.15in} Introduction\\
\noindent 2. \hspace*{.15in} Modularity: definition and some properties\\
\noindent 3. \hspace*{.15in} The modularity of binomial random graphs $\Gnp$\\
\noindent 4. \hspace*{.15in} The modularity of random regular graphs\\
\noindent 5. \hspace*{.15in} Modularity in small world models\\
\noindent 6. \hspace*{.15in}  Modularity and partially observed graphs\\
\noindent 7. \hspace*{.15in}  Modularity and stochastic block models\\ 
\noindent 8. \hspace*{.15in} The modularity of random graphs embeddable in a surface\\
\noindent 9. \hspace*{.15in} Summary of known modularity values\\
References

\vspace*{.1in}

\begin{center}
\section{\hspace{-.3in} . \hspace{.1in} Introduction} \label{sec.intro}
\end{center}
\vspace{-0.15in}
\noindent The ever-increasing amount of network data available in many fields has led to great interest in understanding network structure. In particular we want to be able to identify whether a network can be decomposed into dense clusters or `communities'. 
Modularity was introduced in 2004 by Newman and Girvan~\cite{NewmanGirvan}, to give a measure of how well a graph can be divided into communities, and it now forms the backbone of the most popular algorithms used to cluster real data~\cite{popular}, including the algorithms Louvain~\cite{louvain} and Leiden~\cite{traag2019leiden}. Here a `community' is a collection of vertices which are more densely interconnected than one would expect from looking at vertex--degrees.

There are many applications, including protein discovery, gene co-expression, and identifying connections between websites; see~\cite{fortunato2016community}, \cite{lambiotte2021modularity} and~\cite{porter2009communities} on the use of modularity for community detection in networks.
For example, the paper \cite{elgesem} studied the choice of words used to discuss climate change online. Its authors collated about 3,000 blogs
and, as a preliminary step, produced a graph with a vertex for each blog and an undirected edge for each link between blogs; and then, using the Louvain algorithm~\cite{louvain} (based on modularity), found a partition of these blogs. The authors then analysed how language use varied between these groups.

In the next section we give definitions and notation for modularity, and some basic properties; and we then discuss modularity for several different models of random graphs. 
In Section~3 we discuss the binomial random graph $G_{n,p}$, where edges appear independently with probability $p$, and Section~4 concerns random regular graphs.
In Section 5 we discuss three `small-world' models of random graph: the preferential attachment model, the spatial preferential attachment model, and random graphs in the hyperbolic plane.
In Section 6 we discuss partially observed or `percolated' graphs: mostly we assume that there is a large (unknown) underlying graph $G$, and the edges of $G$ appear independently with some probability $p$, forming the observed graph $G_p$: what can we say about the modularity $\q(G)$?
Section 7 concerns stochastic block models, where there is a `planted partition' which should give a good community structure.  
In Section 8 we consider random planar graphs, and more generally random graphs embeddable in a given surface. 
Finally, Section 9 
gives a summary of known modularity results for different classes of graphs.

\vspace{-.4in}

\begin{center}
\item{}
\section{\hspace{-.3in} . \hspace{.1in}  Modularity: definition and some properties}
\label{sec.mod}
\end{center}
\vspace{-0.15in}
\noindent Given a graph $G=(V,E)$, modularity gives a score to each partition of the vertex-set $V$, and the (maximum) modularity $\q(G)$ of $G$ is the maximum of these scores over all vertex-partitions.
Let~${\mathbf 1}_{vw\in E}$ be the indicator that $vw$ is an edge.
For a set $A$ of vertices, let $e(A)$ (or $e_G(A)$, if needed for clarity) be the number of edges within $A$, and let the volume $\vol(A)$ be the sum of the degrees~$d_v$ over the vertices $v$ in $A$. 
Let $G$ be a graph with $m\geq 1$ edges. Then, for a vertex-partition $\cA$ of~$G$, the \emph{modularity score} of $\cA$ on $G$ is 
\begin{equation}\label{eq.defn} 
 q_\cA(G) =  
\frac{1}{2m}\sum_{A \in \cA} \sum_{v,w \in A} 
\left( {\mathbf 1}_{vw\in E} - \frac{d_v d_w}{2m} \right)
 = 
\frac{1}{m}\sum_{ A \in \cA} e(A) - \frac{1}{4m^2}\sum_{A \in \cA} \vol (A)^2\,.
\end{equation}
The \emph{modularity} of $G$ is $\q(G)=\max_\cA q_\cA(G)$, where the maximum is taken over all vertex-partitions~$\cA$. (For a graph $G$ with no edges, where $m=0$, we set $q_{\cA}(G)=0$ for each vertex-partition~$\cA$, and $\q(G)=0$.) Taking account of the `strength' of connections can be important, and the definition extends naturally to weighted graphs, as we shall see in Section~\ref{sec.edge-sampling} -- note in particular equation~\eqref{eq.defn_w}.

It follows directly from the definition that $0 \leq \q(G) <1$ for all graphs~$G$. Note that isolated vertices are irrelevant -- they are not counted in the formula for the modularity score -- and so we may assume without further comment that our graphs have no isolated vertices.
The second expression for $q_{\cA}(G)$ expresses it as the difference between
the \emph{edge contribution} or \emph{coverage}\, $q^E_\cA(G)= (1/m)\sum_A e(A)$, and the \emph{degree tax}\,  $q^D_\cA(G)= (1/4m^2) \,\sum_A\vol(A)^2$.
By the Cauchy--Schwarz inequality, if $x_1,x_2,\ldots,x_k \geq 0$, then $\sum_i x_i^2 \geq (\sum_i x_i)/k$: thus, if the partition $\cA$ has $k$ parts, then $q^D_\cA(G) \geq 1/k$, and so $q_\cA(G) \leq 1-1/k$.

The original rationale for the definition~\cite{NewmanGirvan} was that while rewarding the partition for capturing edges within the parts, we should penalise by (approximately) the expected number of such edges. We can see this in two different ways, given a suitable degree-sequence ${\bf d}=(d_1, d_2, \ldots, d_n)$. By using the configuration model (see for example~\cite{purplebook} and Chapter~4), we can  generate nearly uniformly at random a multigraph $R$ on $[n]=\{1,2,\ldots,n\}$ with degree-sequence ${\bf d}$ (where multiple edges and loops are allowed); and then the expected number of edges in $R$ between distinct vertices $v$ and $w$ is $d_v d_w/(2m-1)$, where $\sum_v d_v = 2m$.
If instead we sample $G$ uniformly at random from all (simple) graphs with degree-sequence $\mathbf{d}$, then the probability that $vw$ is an edge is very close to~$d_vd_w/2m$ -- see also Theorem~\ref{thm.approxdudv} at the end of this section.  Also, see~\cite{vdense} for a fuller story concerning the degree tax.

We shall discuss the modularity $\q(G)$ of various sorts of random graph $G$. This is very different from considering the modularity of a random partition of a fixed graph.  For let $G$ be any fixed graph (with at least one edge).  Fix $k \geq 2$, and suppose that we generate a random partition $\cal A$ of the vertices by placing each vertex into one of $k$ parts independently with probability $1/k$.  Then, as noted in the concluding remarks of~\cite{treelike},
\begin{equation} \label{eq.randpart}
{\mathbb E}[q_{\cal A}(G)] <  0.
\end{equation}
Thus, comparing the modularity of a partition which we have found to that of a random partition is likely to give a false positive. To see why~\eqref{eq.randpart} holds, observe that, for an edge $e$ in $G$, the probability that both endvertices are placed in the same part is $1/k$, so ${\mathbb E}[q_{\cal A}^E(G)]=1/k$; and the degree tax is always at least $1/k$, and sometimes is larger (see just before Theorem~\ref{prop.1overk} below).

As Reichardt and Bornholdt~\cite{trulymodular} have remarked:
`A differentiation between graphs which are truly modular and those which are not can $\cdots$  only be made if we gain an understanding of the intrinsic modularity of random graphs.'
In the following sections we discuss modularity for several different models of random graphs.

A basic observation is that, given a graph $G$ without isolated vertices, in each optimal partition~$\cA$, each part has size at least $2$~\cite{w1hard} and induces a connected subgraph of $G$.
Also, $V(H)$ is a part of $\cA$ for each connected component $H$ of $G$ with modularity $\q(H)=0$ (see~\cite[Corollary 1.5]{modexpansion}, and see also Theorem 1.4 in that paper for a fuller account).  For example, all complete graphs, `nearly complete' graphs (with $n$ vertices and at least $n(n-2)/2$ edges), and complete multipartite graphs have modularity 0 (see~\cite{vdense} and the references therein).

The idea of a resolution limit was introduced by Fortunato and Barthélemy~\cite{FortBart2008} in 2007. 
Let $G$ be an $m$-edge graph without isolated vertices, which has a connected component $C$ with $e(C) < \sqrt{2m}$: then every optimal partition for $G$ has one part $V(C)$, so all the vertices of $C$ are clustered together.
\vspace{.05in}

\begin{figure}[h!]
    \centering
    \includegraphics[scale=0.7]{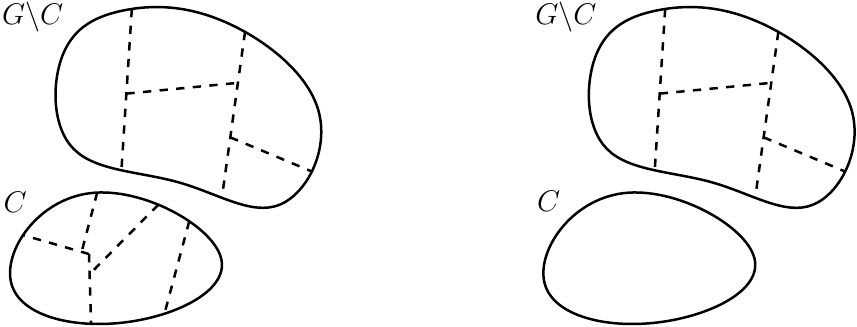}
     \caption{ A partition of $G$ in which $C$ is split into parts (left), and the corresponding partition where~$C$ forms a single part (right).
    }\label{fig.partitionC}
\end{figure}

To see this, consider an optimal partition $\cA$ for $G$, and note that for each part $A$ in $\cA$, the induced subgraph $G[A]$ is connected, and so $A \subseteq V(C)$ or $A \cap V(C)=\emptyset$.
Suppose that there are $k \geq 2$ parts $A_1,A_2,\ldots,A_k$ with $A_i \subseteq V(C)$.
Let $\cB$ be the partition formed from $\cA$ by replacing these parts $A_i$ by their union~$V(C)$ -- see~Fig.~\ref{fig.partitionC}.  Then $\, q_{\cB}^E(G) - q_{\cA}^E(G) \geq (k-1)/m$,
since any connected graph with $k$ vertices has at least $k-1$ edges. 
Also, using again the simple result that if $x_1,x_2,\ldots,x_k \geq 0$ then $\sum_i x_i^2 \geq (\sum_i x_i)/k$, we have
\begin{eqnarray*}
 q_{\cB}^D(G) - q_{\cA}^D(G) 
 &=&
 (\vol(C)^2 - \sum_i \vol(A_i)^2)/ \vol(G)^2 \\
& \leq &
(\vol(C)^2 - \vol(C)^2 /k)/ \vol(G)^2 \; = \; (k- 1)/k \;\, e(C)^2/m^2.
\end{eqnarray*}
Hence
\[
 m\,(q_{\cB}(G) - q_{\cA}(G))  \geq (k-1) (1- e(C)^2/km) > 0 \;\;\; \mbox{ since } \: e(C) < \sqrt{km}.
\]
Thus, no optimal partition splits $C$, and we have proved the resolution limit result.

Consider an example where the component $C$ consists of any two connected graphs $C_1$,~$C_2$ with $e(C_1)=e(C_2)$, joined by a single edge. As we have seen, if $\sqrt{2m} > e(C)$ then $C$ is never split in an optimal partition. But if $\sqrt{2m} < e(C)$, then $C$ is always split: to see this, note that if $\cA$ splits $V(C)$ into $V(C_1)$ and $V(C_2)$, then in the last displayed equation, the first inequality holds at equality.
This example also highlights the sensitivity of optimal partitions to noise in the network: if that edge between $C_1$ and $C_2$ (perhaps a mistake in the data) had not been there, then $C$ would be split in every optimal partition. In contrast, although the structure of optimal partitions is not robust to small changes in the set of edges, the modularity value of the graph is robust in this sense, as we see in the following result  (see~\cite{sampling}, where it is also shown that the constant 2 on the right side of the inequalities is best possible). 
\begin{thm}
\label{lem.robustness}
Let the graph $G=(V,E)$ have $m \geq 1$ edges, and let $\cA$ be a partition of $V$. If $E_0$ is a non-empty subset of $E$, and $G'=(V, E \setminus E_0)$, then  
\[ |q_\cA(G)- q_\cA(G')| < \frac{ 2 \, |E_0|}{m} \;\; \mbox{ and } \;
\; |\q(G)-\q(G')|< \frac{2 \, |E_0|}{m}.\]
\end{thm}
\vspace{-0.05in}
\noindent
Next let us consider the number of parts in the partitions.  
We noted that $q_\cA(G) \leq 1-1/k$ for any $k$-part partition $\cA$.
The following result of Dinh and Thai~\cite{dinh2011finding} shows that if $\q(G)>0$, then there is a 2-part partition (a bipartition) $\cB$ for $G$ with $q_\cB(G)>0$. 
\begin{thm}\label{prop.1overk}
Let $G=(V,E)$ be a graph, and let $k \geq 2$. Then there is a partition $\cB$ with at most~$k$ parts such that $q_\cB(G) \geq (1-\frac1{k}) \, \q(G)$. 
\end{thm}
\vspace{-0.1in}
\noindent
Theorem~\ref{prop.1overk} is best possible for each $k \geq 2$, in that we cannot replace the factor $(1-\frac1{k})$ by any larger value. For, $q_\cB(G) \leq 1-\frac1{k}$ since $\cB$ has at most $k$ parts; and, given $\eps>0$, there is a graph $G$ such that $(1-\frac1{k} + \eps)\, \q(G) > 1-\frac1{k}$.  (For example, we may consider $G$ as consisting of many disjoint edges.)

The deterministic result Theorem~\ref{prop.1overk} has a simple and beautiful probabilistic proof. Let $\cA=(A_1,A_2,\ldots,A_t)$ be a $t$-part partition of $V$ with $q_\cA(G)=\q(G)$, where we may assume that $k<t$.
Let the random variables $X_1,X_2,\ldots,X_t$ be independent, with each uniformly distributed on $[k]$. We obtain a random partition $\cB=(B_1,B_2,\ldots,B_k)$ of $V$ (with at most $k$ non-empty parts) by setting $B_i = \bigcup_{\{j:X_j=i\}} A_j$.

Let $G$ have $m$ edges; and for $v,w \in V$ let $C_{v,w} = {\bf 1}_{vw \in E} - d_vd_w/2m$, and (in a mild abuse of notation) write $vw \in \cA$ if $v$ and $w$ are in the same part of $\cA$. (Note that we allow $v=w$ here.) Now
\[ q_\cA(G) = \tfrac1{2m} \sum_{v,w} 
{\bf 1}_{vw \in \cA}  \, C_{v,w}
 \;\; \mbox{ and } \;\; \E[q_{\cB}(G)]  = 
  \tfrac1{2m}\,  \sum_{v,w} \E[
{\bf 1}_{vw \in \cB}] \, C_{v,w}. \]
Note that $\,\E[{\bf 1}_{vw \in \cB}] =  {\bf 1}_{vw \in \cA} + \tfrac1{k}\, {\bf 1}_{vw \not\in \cA}$.
Also, $\sum_{v,w} C_{v,w} = 2m- (2m)^2/2m = 0$, and so
\begin{equation} \label{eq.sumC} \nonumber
 \sum_{v,w} {\bf 1}_{v w \not\in \cA}\, C_{v,w} = - \sum_{v,w} {\bf 1}_{vw \in \cA}\, C_{v,w}\,.   
\end{equation}
Hence
\[ \E[q_{\cB}(G)] \; = \;
\tfrac1{2m} \, \sum_{v,w} \big({\bf 1}_{vw \in \cA} + \tfrac1{k} {\bf 1}_{vw \not\in \cA} \big) \, C_{v,w} \; = \;
(1-\tfrac1{k})\, q_{\cA}(G) \; = \; (1-\tfrac1{k})\, \q(G) \,,\]
and this yields Theorem~\ref{prop.1overk}.

\smallskip

Again let $G=(V,E)$ be a graph with $m \,(\geq 1)$ edges.
The \emph{relative modularity} of a non-empty set $A \subseteq V$ is
\[ q_{r}(A) = \frac{2e(A)}{\vol(A)} - \frac{\vol(A)}{2m}.\]
Then, for any partition $\cA$ of $V$, $q_\cA(G)$ is a weighted average of values $q_r(A)$,
\begin{equation}\label{eq.rel_mod} q_\cA(G) = \sum_{A \in \cA} \frac{\vol(A)}{2m} q_r(A) .\end{equation}
Denote the maximum value of $q_r(A)$ 
over non-empty $A \subseteq V$ by $\q_r(G)$.  Then, by~\eqref{eq.rel_mod},
$q_\cA(G) \leq \q_r(G)$, and so $\q(G) \leq \q_r(G)$. This inequality often yields our best bounds on $\q(G)$  -- see for example the proof of Theorem~\ref{thm.cubic}. 
In contrast, however, we can have $\q(G) =o(1)$ and $q^*_r(G) = 1-o(1)$. For example, consider an $n$-vertex graph $G$ consisting of about $\log n$ disjoint edges and a complete graph on the remaining vertices. 
Then the number of edges `missing' from $G$ is $o(n^2)$ and $\q(K_n)=0$ (see~\cite{nphard}), and so $\q(G)=o(1)$ by Theorem~\ref{lem.robustness}. However, $q_r(A)=1-o(1)$ when $A$ is the set of vertices of the disjoint edges.

The following result concerning simple random graphs on $[n]$ (which we mentioned earlier) is due to Liebenau and Wormald~\cite[Theorem~1.6]{liebenau2024asymptotic} (and 
more precise results are given in Lemma~7.1 of the same paper).
\begin{thm}\label{thm.approxdudv}
Let $m=m(n) \in \N$, and assume that $d=2m/n$ satisfies $(\log n)^\omega < d < \mu_0 n$, where $\mu_0>0$ is sufficiently small and $\omega=\omega(n) \to \infty$ as $n \to \infty$. Let $\alpha < 3/5$, and let the $n$-vector $\textbf{d} \in \N^n$ with $\sum_i d_i=2m$ satisfy the `near-regular' condition $\max_i |d_i - d| \leq d^\alpha$. Now suppose that 
we sample~$G$ uniformly from all graphs with degree-sequence~$\textbf{d}$. Then, for distinct vertices $v$ and~$w$,
\[
\Pr(vw \in E(G)) = \frac{d_v d_w}{d(n-1)}  \Big( 1 - \frac{(d_v - d)(d_w - d)}{d(n-1-d)}\Big) + O\bigg( \frac{\sqrt{d \log n}}{n^2} + \frac{(\sqrt{d \log n})^3}{n^3}
\bigg).\]
\end{thm}

\noindent Noting that $d(n-1)=2m(1-d/2m)=2m(1-1/n)$ 
this result implies 
\[
\Pr(vw \in E(G)) = \frac{d_v d_w}{2m}   \Big( 1 - \frac{(d_v - d)(d_w - d)}{d(n-1-d)}\Big)\Big(1-\frac{1}{n}\Big)^{-1}  
+ O\Big( \Big(\frac{\log n}{n} \Big)^{3/2} \Big).
\]
Since by assumption $|(d_v - d)(d_w - d)|/d\leq d^{2\alpha -1} = o(n^{1/5})$,
we obtain the simpler equation
\[
\Pr(vw \in E(G)) = 
\frac{d_v d_w}{2m}\Big( 1+ o\Big(\frac{1}{n^{4/5}}\Big)\Big) 
+ O\Big( \Big(\frac{\log n}{n} \Big)^{3/2} \Big).\]
\vspace{-0.1in}
\noindent

\vspace*{.2in}
\begin{center}
\section{\hspace{-.3in} . \hspace{.1in} The modularity of binomial random graphs {$\boldsymbol{G_{n,p}}$} 
}
\label{sec.binom}
\end{center}
\vspace{-0.15in}
\noindent Let~$n$ be a positive integer.  Given $0 \leq p \leq 1$, the \emph{binomial random graph} $\Gnp$ has vertex-set~$[n]$ and the $\binom{n}{2}$ possible edges appear independently with probability $p$.  Given an integer~$m$ with $0 \leq m \leq \binom{n}{2}$, the \emph{Erd\H{o}s--R\'enyi random graph} $\Gnm$ is sampled uniformly from the graphs on vertex-set $[n]$ with $m$ edges.  These two random graphs are closely related when $m$ is about $\binom{n}{2} p$. We shall focus here on~$\Gnp$, and mention $G_{n,m}$ only briefly towards the end of this section.

Our first theorem, the `three phases theorem', is from~\cite{ERmod} and gives the big picture. The three phases correspond to when the expected vertex-degree (essentially $np$) is (a) at most about 1, (b) larger than 1 but bounded, or (c) tending to infinity.

\needspace{8\baselineskip}
\begin{thm}\label{thm.usER}
Let $p=p(n)$ satisfy  $0 < p \leq 1$. 

\noindent {\rm (}a{\rm )} If $n^2p \rightarrow \infty$ and $np\leq 1+o(1)$, then whp $\; \q(\Gnp) = 1-o(1)$.

\noindent {\rm (}b{\rm )} Given constants $1< c_0 \leq c_1$,  there exists $\delta=\delta(c_0,c_1)>0$ such that, if $c_0 \leq np \leq c_1$ for $n$ sufficiently large, then whp $\; \delta<\q(\Gnp)<1-\delta$. 

\noindent {\rm (}c{\rm )} If $np\rightarrow \infty$ then whp $\; \q(\Gnp)=o(1)$.
\end{thm}
\vspace{-0.1in}
\noindent
For any graph $G$, let $\cC$ be the partition of the vertex-set into the vertex-sets of its connected components; then $q_\cC^E(G)=1$. If $np$ is a little less than 1, then whp $\cC$ is the unique optimal partition for $\Gnp$ (see~\cite{thesis}).
Also, if $np \leq 1+o(1)$, then $q_\cC(\Gnp) \overset{p}\rightarrow 1$, and so $\q(\Gnp) \overset{p}\rightarrow 1$.  To show this, it suffices to check that whp the fraction of edges in any component is $o(1)$, and so the degree tax is~$o(1)$; this yields part (a) of Theorem~\ref{thm.usER}. Parts (b) and (c) largely follow from more detailed results as in Theorem~\ref{thm.growthRate} below, see~\cite{ERmod} for the upper bound in part (b). (We shall see later that part (c) also follows from Theorem~\ref{thm.moddiff}\,(a) on percolated random graphs.)

In each of the parts (a), (b), (c) of Theorem~\ref{thm.usER} we can be more precise. For the rest of this section, we shall focus on the denser ranges (b) and (c). The following result is from~\cite{ERmod}. 
\begin{thm}\label{thm.growthRate}
There is a constant $b$ such that, for all  $0< p=p(n) \leq 1$, we have $\q(\Gnp)< b/\sqrt{np}$ whp. 
Also, given $0<\eps<1$, there exists $a =a(\eps) >0$ such that, if $p=p(n)$ satisfies $np\geq 1$ and $p \leq 1-\eps$ for $n$ sufficiently large, then  $\q(\Gnp) > a/\sqrt{np}$ whp. 
\end{thm}
\vspace{-0.1in}
\noindent
This theorem shows that, if $p=p(n)$ satisfies $1/n \leq p \leq 0.99$, then $\q(\Gnp) = \Theta((np)^{-1/2})$ whp, confirming a conjecture in the physics literature~\cite{trulymodular} (see also the comment following Theorem~\ref{thm.lowerSqrt2} below).

The lower bound in Theorem~\ref{thm.growthRate} comes from analysing a simple algorithm `Swap', which works as follows. Given a graph $G$ with vertex-set $[n]$, we start with the bipartition $(A,B)$ into the odd-numbered and even-numbered vertices: this partition has modularity score near 0 whp.  By interchanging some pairs $a_i, b_i$ between $A$ and $B$, whp we can significantly increase the edge contribution, without changing the distribution of the degree tax.
In a little more detail, assume for simplicity that $n$ is divisible by 6, 
let $k=n/6$, and relabel $A$ as $a_1,a_2,\ldots,a_{3k}$ and $B$ as $b_1,b_2, \ldots,b_{3k}$. 
We keep $A_1=\{a_{2k+1}, a_{2k+2},\ldots,a_{3k}\}$ and $B_1=\{b_{2k+1}, b_{2k+2},\ldots,b_{3k}\}$ fixed.
For each $i=1,2,\ldots,2k$ we consider the edges between $a_i,b_i$ and  $A_1, B_1$, and interchange (swap) $a_i$ and~$b_i$ if that would increase the number of these edges which lie within the parts.
See also Section~\ref{sec.stoch-block} below where we discuss the stochastic block model.

Now consider the upper bound in Theorem~\ref{thm.growthRate}, which will follow from spectral results. Define the \emph{normalised Laplacian} of a graph $H$ to be $\mathcal{L}_{H}= \matr{I} - \matr{D}^{-1/2}\matr{A} \matr{D}^{-1/2}$, where $\matr{A}$ is the adjacency matrix of $H$, $d_i$ is the degree of vertex $i$, and $\matr{D}^{-1/2}$ is the diagonal matrix with $i$th diagonal entry~$d_i^{-1/2}$.  The \emph{spectral gap} of an $n$-vertex graph $H$ is $\bar{\lambda}_{H} = \max_{i \neq 0} |1-\lambda_i|$, where $0=\lambda_0 \leq \lambda_1 \leq \cdots \leq \lambda_{n-1}$ are the eigenvalues of~$\mathcal{L}_{H}$. 

The modularity of a graph $G$ is bounded above by its spectral gap: we have $\q(G) \leq \bar{\lambda}_G$ (see~\cite{ERmod}). Moreover,~$\q(G)$ is bounded above by  a simple function of the relative size and the spectral gap of any subgraph $H$, as in the following result~\cite[Corollary 1.9]{modexpansion}. 

\begin{thm}\label{thm.spectralub}
Let the graph~$G$ have  a subgraph~$H$ with spectral gap $\bar{\lambda}_H$, and let $\alpha = e(H) / e(G)$. Then
\begin{equation*} q^*(G) \leq 1- \alpha \min \{\alpha, 1-\bar{\lambda}_H \} .\end{equation*}
\end{thm}
\vspace{-0.1in}
\noindent 
To complete the proof of the upper bound on $\q(G_{n,p})$, we apply the bound in~Theorem~\ref{thm.spectralub} to a graph $H$ obtained from $\Gnp$ by deleting a few vertices, such that whp $H$ has small spectral gap.

\medskip

In 2025, Rybarczyk and Sulkowska~\cite{rybarczyksulkowska2025Gnp} improved the upper bound in Theorem~\ref{thm.growthRate}, by replacing the constant $b$ by the explicit value 3.06, showing that
\begin{equation} \label{eqn.b}
    \mbox{ for all } 0< p=p(n) = o(1), \;\; \q(\Gnp)< 3.06/\sqrt{np} \;\; \mbox{whp}. 
\end{equation}
In their proof they first used Chernoff's inequality to show bounds on $e(S)$, $e(\bar{S})$ and $e(S, \bar{S})$ which whp hold simultaneously for all vertex subsets $S$.

If we weaken the lower bound in Theorem~\ref{thm.growthRate} by introducing an extra factor $\sqrt{1\!-\!p}$, then we can extend the bound for $p$ up to when the complementary graph has expected average degree just above 1 (see~\cite{ERmod}).
Observe that if $p$ is bounded below 1 (as in the lower bound in Theorem~\ref{thm.growthRate}), then this factor is just a constant. The next two results are from~\cite{vdense}.
\begin{thm} \label{thm.from-modER}
For each $\eps>0$, there exists $\alpha>0$ such that, if $p=p(n)$ satisfies $1/n \leq p \leq 1-(1+\eps)/n$, then $\q(G_{n,p}) \geq \alpha \sqrt{(1\!-\!p)/np}$ whp.
\end{thm}
\vspace{-0.1in}
\noindent
We saw in Theorem~\ref{thm.growthRate} that $\q(\Gnp)< b/\sqrt{np}$ whp, and so if $np \to \infty$ then $\q(\Gnp) \to 0$ whp (as in Theorem~\ref{thm.usER}(c)).
When do we have the stronger result that whp $\q(\Gnp)=0$?  
Consider the very dense case, where $p=1- \Theta(1/n)$.  
In this case, 
$\sqrt{{(1-p)}/{np}}= \Theta(1/n)$, and so the lower bound in part (c) below follows directly from Theorem~\ref{thm.from-modER}. 
\begin{thm}\label{thm.Gnp}
Let $0 \leq p=p(n)\leq 1$. 
\\ 
\rm{(}a\rm{)}
If $p$ satisfies $p \geq 1- 1/n + \omega(n)\, n^{-3/2}$,
where $\omega(n) \to \infty$ as $n \to \infty$, then $\q(G_{n,p})=0$ whp.\\
\rm{(}b\rm{)} If $p=1-1/n$, then $\mathbb{P}\big(\q(G_{n,p})=0\big)\rightarrow 1/2$ as $n \to \infty$.\\
\rm{(}c\rm{)}
Given $1<c_1<c_2$, there exist $\alpha, \beta$ with $0< \alpha < \beta$ such that, if $p$ satisfies 
$1- c_2/n \leq p \leq 1-c_1/n$,
then $\alpha/n \leq \q(G_{n,p}) \leq \beta/n$ whp.
\end{thm}
\vspace{-0.1in}
\noindent
 It is known that (for $n \geq 4$) the least number of edges missing from an $n$-vertex graph $G$ with $\q(G)>0$ is $\lfloor n/2 \rfloor+1$ (see~\cite{vdense}). 
 But when $1-p$ is very small (below $1/n$), then whp at most~$n/2$ edges are missing from $\Gnp$ and so
 $\q(G_{n,p})=0$; this yields part (a) of Theorem~\ref{thm.Gnp}. Part~(c) says that when~$1-p$ is just a little larger (above $1/n$), then whp $\q(G_{n,p}) = \Theta(1/n)$.

\needspace{3\baselineskip}
\begin{center}
 {\bf The Erd\H{o}s--R\'enyi random graph $\Gnm$}
\end{center}

\noindent Recall that, given $m=m(n)$, the random graph $\Gnm$ is sampled uniformly from the graphs on~$[n]$ with $m$ edges.  It is similar to $\Gnp$ with $p=m/\binom{n}{2}$, and so it is natural that there are results on $\q(\Gnm)$ which correspond to Theorems~\ref{thm.usER} to~\ref{thm.Gnp} above (see~\cite{ERmod}, Section~1.3). For example, if $m \to \infty$ and $m \leq (\frac12 + o(1))n$, then $\q(G_{n,m}) \to 1$ in probability, and if $m/n \to \infty$, then $\q(G_{n,m}) \to 0$ in probability.

\needspace{3\baselineskip}
\begin{center}
 {\bf Friendly bisections}
\end{center}

\noindent Given a fixed $p$ with $0<p<1$, whp $\q(G_{n,p}) = \Theta(n^{-1/2})$ by Theorem~\ref{thm.growthRate}, and so whp the maximum score over bipartitions is $\max_{|\cA|=2} q_{\cA}(G_{n,p}) = \Theta(n^{-1/2})$ by Theorem~\ref{prop.1overk}. Thus, whp there is a bipartition with small but positive modularity score.

A different take on bipartitions may reflect the community structure of the graph. 
A partition is \emph{balanced} if the sizes of the parts differ by at most 1, and a balanced bipartition is also called a \emph{bisection}.
A \emph{friendly bisection} of an $n$-vertex graph is a  bisection in which all but $o(n)$ vertices have more neighbours in their own part than in the other part.  It is shown in~\cite{ferber2022friendly} that whp $G_{n,1/2}$ has a friendly bisection.

A bisection such that all but one vertex in each part have more neighbours in their own part may have modularity score 0. For example, let the graph $G$ consist of two disjoint wheels on $k+1$  vertices, so that $G$ has $n=2k+2$ vertices and $m=4k$ edges.
Let $A$ contain the rim of one wheel and the centre of the other, let $B$ be the complementary set of vertices, and consider the bisection $\cA=\{A,B\}$.
Each of the $2k$ vertices in the rims has two neighbours in the same part and one in the other part, and the two centres have all $k$ neighbours in the other part.  Thus, $\cA$ is a friendly bisection of $G$; but $q_\cA^E(G)= 2k/m = \frac12$, and $q_\cA^D(G) =\frac12$ by symmetry, and so $q_\cA(G)=0$.

However, suppose that, given a bipartition $\cA$ for $G$ (not necessarily balanced), \emph{all} vertices have at least as many neighbours in their own part as in the other part, 
and at least one vertex has strictly more. In this case,
$q_\cA(G)>0$. To see this, let there be $x$ edges between $A$ and $B$, $y$ edges within $A$, and $z$ edges within $B$. 
Then $x \leq 2y$ and $x \leq 2z$, with at least one of these inequalities strict, and so $x^2< 4yz$.
But
$q_{\cA}^E(G) = (y+z)/(x+y+z)$ and
$q_{\cA}^D(G) = \big( (x+2y)^2 + (x+2z)^2 \big)/(2x+2y+2z)^2$.  Hence 
\[
4(x+y+z)^2\, q_{\cA}(G) =
4(x+y+z)(y+z) - \big( (x+2y)^2 + (x+2z)^2 \big) = 2(4yz - x^2) >0\,, \]
and so indeed $q_{\cA}(G)>0$.
\medskip

\begin{center}
\section{\hspace{-.3in} . \hspace{.1in} The modularity of random regular graphs}
\label{sec.reg}
\end{center}
\vspace{-.15in}
For a regular graph $G$, modularity scores have a simplified form (see~\eqref{eq.defn} and \eqref{eq.rel_mod}): 
\begin{equation}
    \label{eq.regmod}
 q_\cA(G) = \sum_{A \in \cA} \big( \frac{e(A)}{m} - \frac{|A|^2}{n^2} \big)  = \sum_{A \in \cA} \frac{|A|}{n}\,q_r(A),\, \;\; \mbox{ where } \;\; q_r(A) = \;\big( \frac{2 e(A)}{r|A|} - \frac{|A|}{n} \big).
\end{equation}
Random 2-regular graphs $G_{n,2}$ are easy to handle. They have only a few cycles whp, and so, by the robustness result in Theorem~\ref{lem.robustness}, whp their modularity is close to that of the $n$-cycle $C_n$.  It follows (see~\cite{treelike}) that
\begin{equation}\label{eq.2_regular} \q(G_{n,2}) = 1 - \frac2{\sqrt{n}} +o\big(\frac{\log^2n}{n}\big)
\;\;\; \mbox{ whp.}
\end{equation}
In~\cite{treelike} it is also shown that, for random 3-regular graphs $G_{n,3}$ (with $n$ even),
\[  2/3 - o(1)\leq \q(G_{n,3}) < 0.804 \;\;\; \mbox{ whp.} \]
Here, the lower bound follows easily from the result that whp $G_{n,3}$ has a Hamiltonian cycle (see~\cite{robinson1992almost}).
We ignore the rest of the graph (a perfect matching) and cut the circuit into $\omega(n)$ paths of about equal length to give a vertex-partition, with edge contribution $\frac23 -o(1)$ and degree tax~$o(1)$.
Indeed, for each $r \geq 2$, 
$\q(G) \geq 2/r - 2 \sqrt{6/n}$ for {\em every} $n$-vertex $r$-regular graph $G$ (see~\cite[Proposition 1.4]{treelike}.)
It was conjectured in~\cite{treelike} that whp $\q(G_{n,3}) > 2/3 + \delta$ for some $\delta>0$.
This was shown in~\cite{lichev2022modularity} together with an improved upper bound. 
\begin{thm} \label{thm.cubic}
\;\;\;\; $0.667026 \leq \q(G_{n,3}) < 0.789998$ \mbox{ whp}.
\end{thm}
\vspace{-0.1in}
\noindent
To prove that $\q> \tfrac23 +\delta$ it suffices to show that there exists a set $C$ of $\Theta(n)$ vertices such that  
$q_r(C) \geq \frac23+\eps$ 
and the remaining graph may still be partitioned so that the weighted sum contributes $(\frac23 - o(1))(1-|C|/n)$ to the modularity score (see~\eqref{eq.regmod}). We noted earlier that to prove an upper bound, it suffices to show that the relative modularity of all vertex-sets has at most some value.

To construct a vertex-set with high relative modularity and linear size, and thus to prove the lower bound in Theorem~\ref{thm.cubic}, we consider the configuration model. 
Choose a vertex $v$ uniformly at random, and start exploring its connected component one edge at a time in a certain order, until a set $C_0$ of $\eps n$ vertices has been explored. 
Next, explore one by one the open half-edges sticking out of $C_0$, to grow it little by little. Then start adding to $C_0$ short paths containing exactly two explored vertices -- the first and the last vertex of each path. In this way, the relative modularity of~$C_0$ increases at each step. By analysing the first few steps of this procedure, using the differential equation method~\cite{cookbook} carefully, and optimizing over $\eps$, we can complete the proof.

To prove the upper bound in Theorem~\ref{thm.cubic}, we note that, for the modularity to be at least the given upper bound, there must be a set $A$ of vertices having relative modularity at least that value. 
First we show that there is an $\eps_0>0$ such that whp this does not hold for any `small' set $A$ with $|A| < \eps_0 n$.
Then we show that it suffices to prove that there is no `large' set $B$, containing many edges, and such that the induced subgraphs on $B$ and on its complement contain only vertices of degree~2 or~3. Finally, we use the first moment method to show that whp there is no such set $B$.

The next theorem (from~\cite{treelike}) shows that $\q(G_{n,r})$ is $\Theta(1/\sqrt{r})$ whp for large (fixed) values~$r$. For simplicity, we restrict our attention to even values of $n$.

\begin{thm}\label{thm.rreg_larger}
For each $r \geq 2$,  \,$\q(G_{n,r}) < 2\, / \sqrt{r} \, \mbox{ whp}$,
and there is a constant $r_0$ such that, for each $r\geq r_0$, \,$\q(G_{n,r}) > 0.7632 \, / \sqrt{r} \, \mbox{ whp}$.
\end{thm}
\vspace{-0.1in}
\noindent
Here, the lower bound follows from a result on the bisection width of random regular graphs. Dembo, Montanari and Sen~\cite{demboMontanariSen} showed that whp there is a vertex-set $U$ with $|U|=n/2$ such that
$e(U, \bar{U}) \leq rn/4 - c \sqrt{r}n/2+o(\sqrt{r}n)$,
where $c=0.76321\,\pm\,0.00003$.
Thus, in this case, the bipartition $\cA=\{U,\, \bar{U}\}$ has edge contribution at least 
$q^E_\cA= \tfrac12 +  (c+o(1))/\sqrt{r}$ and degree tax $q^D_\cA= \tfrac12$, and this establishes the whp lower bound. The upper bound follows from a spectral bound on random regular graphs by Friedman~\cite{friedman2008proof}, and from the inequality in Theorem~\ref{thm.spectralub} (with $H=G$, $\alpha=1$) which provides an upper bound on the modularity of a graph in terms of its spectral gap. 
\medskip \vspace{-0.1in}

\begin{center}
\section{\hspace{-.3in} . \hspace{.05in} Modularity in small-world models}
\label{sec.smallworld}
\end{center}
\vspace{-0.15in}
In this section we discuss modularity for some `small-world' models, where typical distances are small.
\medskip

\begin{center}
  {\bf The preferential attachment model}
\end{center}
\vspace{-0.1in}
The preferential attachment random graph model $G_n^h$ is of particular interest because it exhibits a `rich-get-richer' behaviour: vertices which already have many incident edges are more likely to be joined to the incoming vertices;
also, vertex-degrees follow a power-law distribution. To construct~$G_n^h$ on the vertex-set $[n]$, we first construct a graph (with loops) on the `mini-vertices'~$1,2,\ldots,hn$. Let $T_1$ be the graph consisting of the single mini-vertex $1$ with a loop $e_1$ (so that $v_1$ has degree~2). For $1 \leq t \leq hn-1$, the graph $T_{t+1}$ is constructed from $T_{t}$ by adding the mini-vertex $t+1$ and joining it by an edge~$e_{t+1}$ to a mini-vertex $s\in \{1,2,\ldots,t\}$ with probability~$d^{t}_s/(2t+1)$, and adding a loop at the mini-vertex $t+1$ with probability $1/(2t+1)$, 
where $d_s^t$ is the degree of the mini-vertex~$s$ in~$T_{t}$. 
Finally, for $h\geq 2$, the random multigraph $G_n^h$ is obtained from $T_{hn}$ by merging each set of~$h$ consecutive mini-vertices together into a single vertex, keeping loops and multiple edges. Observe that $G_n^h$ has~$hn$ edges, and so has average degree $2h$. We consider whp bounds on $ \q(G^h_n)$ for large values of $h$.
\begin{thm}\label{thm.PA}
If $h\geq 2$, then whp  \, $\Omega(1/\sqrt{h}\;) = \q(G^h_n)
    \, = \,  O(\sqrt{\log h}\;  / \sqrt{h} \; )$. 
\end{thm}
\vspace{-0.1in}
\noindent
Here, the lower bound was proved by Prokhorenkova, Pra{\l}at and Raigorodskii~\cite{prokhorenkova2017modularity}, by analysing an algorithm which, as the random graph was constructed, iteratively placed the new vertex into one of two parts, always choosing the part to which that vertex had more incident vertices. 
The same lower bound was also shown to hold for a generalised model of preferential attachment  in~\cite{mod2023universal} -- using only properties of the degree-sequence of that random graph. 
Note that the whp lower bounds for both the preferential attachment model and the binomial (or Erd\H{o}s--R\'enyi) model are~$\Theta(1/\sqrt{d})$, where~$d$ is the expected average degree. 

The upper bound in Theorem~\ref{thm.PA} was proved by Rybarczyk and Sulkowska~\cite{rybarczyksulkowska2025PA}, and we describe their proof. 
The main step is to show that whp, for all vertex-sets $S \subseteq V(G)$, both $e(S)$ and $\vol(S)$ are bounded by functions of a deterministic value $\mu(\tilde{S})$, where $\tilde{S}$ is the set of mini-vertices in $T_{hn}$ corresponding to $S$. Note that the index of a mini-vertex gives the time at which it was added in the construction process. 
If we define 
\[\mu(\tilde{S}) = \tfrac{1}{2}\sqrt{\pi} \sum_{j\in \tilde{S}} (2j-3)!!/(2j-2)!!\,,\] 
we may establish that whp, for all non-empty $S \subseteq V(G)$,
\[ e(S) \leq \mu(\tilde{S})^2 + O\big( n \sqrt{h \log h}\; \big)  \;\;\mbox{ and } \;\; \vol(S) \geq 2\sqrt{hn}\; \mu(\tilde{S}) + O\big(n \sqrt{h \log h}\;\big). \]
We also use the nice observation that
\[ \q(G)\leq 4\max_{S\subseteq V(G)} \left(e(S)/e(G) - \vol(S)^2 / \vol(G)^2 \right)\;\!,\]
which holds since by Theorem~\ref{prop.1overk} there is a bipartition $\cA$ with $\q(G) \leq 2 q_\cA(G)$.
From these inequalities, and bounds on the deterministic function $\mu(\tilde{S})$, the upper bound on $\q(G_n^h)$ quickly follows.

\begin{center}
  {\bf The spatial preferential attachment model}
\end{center}
\vspace{-0.1in}
The spatial preferential attachment model (SPA) yields a sequence of $t$-vertex directed graphs $(G_t)_{0 \leq t \leq n}$, where each vertex has an embedding in $[0,1]^m$. (The edges will be undirected later.) The vertices are added one by one. 
For a vertex $v$ already present in the graph, an edge $vw$ may be added from $v$   
to the newly added vertex $w$ if $w$ falls within a certain distance of $v$, such that existing vertices of higher degree may have an edge to
$w$ at greater distances; here, distances are taken as on a torus. We follow the notation of~\cite{prokhorenkova2017modularity}. Let $a_1,a_2>0$ and let $0<p \leq 1$. The `sphere of influence' $S(v,t)$ of vertex $v$ at time $t$ is defined to be the ball centred at $v$ with volume \[ |S(v, t)| = \min \{ (a_1 d_{G_t}^-(v)+a_2)/{t}, 1 \}, \] where $d_{G_t}^-(v)$ is
the in-degree of $v$ in $G_t$.  Let $G_0$ be the empty graph with no vertices. We construct the graph~$G_t$ with vertex-set $V_t$ from $G_{t-1}$ by adding a new vertex $w_t$ randomly embedded uniformly in~$[0,1]^m$; then, independently for each vertex $v \in V_{t-1}$, if $w_t \in S(v, t-1)$ we 
add the directed edge~$(w_t, v)$ with probability~$p$. 
Finally, let $\widetilde{G}_n$ be the undirected graph constructed from $G_n$ by replacing each directed edge by an edge. Prokhorenkova,  Pra{\l}at and Raigorodskii~\cite{prokhorenkova2017modularity} proved that~$\widetilde{G}_n$ has asymptotically maximal modularity.
\begin{thm} \label{thm.SPA}
Let $m \in {\mathbb N}$, let $a_1, a_2 >0$, let $p\in (0,1]$, and suppose that 
$p\,a_1 < 1$. Then for~$\widetilde{G}_n$ generated as above, \,whp\, $\q(\widetilde{G}_n) = 1 - o(1)$.
\end{thm}
\vspace{-0.1in}
\noindent
The proof of Theorem~\ref{thm.SPA} proceeds by giving an explicit construction for partitioning $[0,1]^m$ into a slowly growing number of slices.
Interestingly, this partition depends on parameters of the model, but is independent of the random graph which has been sampled. The partition of $[0,1]^m$ induces a vertex-partition of $\widetilde{G}_n$, which is then shown to yield modularity score near~1. In the proofs that the edge contribution is near~1 and that the degree tax is near~0, the vertices are considered separately based on their arrival time to the graph.

\needspace{3\baselineskip}
\begin{center}
{\bf The modularity of random graphs in the hyperbolic plane}
\end{center}
\vspace{-0.1in}
We next consider the KPKBV model introduced by Krioukov {\em et al.}~\cite{krioukov2010hyperbolic}, and denote the $n$-vertex graph by $G_{n, \alpha, \nu}$ for parameters $\alpha, \nu>0$. This graph is of interest because it has many properties associated with observed networks: power-law degree distribution with exponent $2\alpha+1$, non-vanishing clustering coefficient and short typical distances between vertices -- see, for example, \cite{fountoulakis2021clustering} and the references therein. Note that an exponent between $2$ and $3$ for the degree distribution corresponds to $\alpha\in (\tfrac12, 1)$, and this is the range most usually considered.

For $\alpha, \nu>0$, we may construct the $n$-vertex hyperbolic random graph $G'_{n, \alpha, \nu}$ as follows.
We generate the position $(r_i, \Theta_i)$ (in polar co-ordinates) of vertex $v_i$, independently for each $i$. 
Set $R=R(n, \nu)$ to satisfy $n=\nu e^{R/2}$,
and sample $r_i$ with probability density function
$\rho_n(r) =\alpha \sinh (\alpha r)/(\cosh (\alpha R ) - 1)$
for $0<r <R$.
Also, sample $\Theta_i$ uniformly in $(0, 2\pi]$, independently of~$r_i$.  Include the edge $v_i v_j$ if the vertices are within (hyperbolic) distance $R$ of each other -- that is, 
if 
\[
\cosh(r_i)\cosh(r_j)-\sinh(r_i)\sinh(r_j) \cos(|\Theta_i - \Theta_j|_{2\pi}) \leq \cosh(R),
\]
where $|a|_b= \min(|a|, b-|a|)$. Define also the  `Poissonisation' $G_{n, \alpha, \nu}$ of the above model -- which has~${\rm Po}(n)$ points, rather than $n$. The vertex-set is the point-set of a Poisson point process on the disc of radius $R$ with intensity $n(2\pi)^{-1}\rho_n(r)dr d\Theta$. The edges are included as for the original model. 

The Poissonisation of the hyperbolic graph has asymptotically maximal modularity value (see~\cite{chellig2022modularity}).  
\begin{thm}
\label{thm.hyperbolic} 
For any $\alpha > \tfrac12$ and $\nu >0$, \,whp\, $\q (G_{n, \alpha, \nu}) = 1-o(1).$
\end{thm}
\vspace{-0.1in}
\noindent
Loosely, this result follows by cutting up the space in a natural way, regardless of the graph, and then showing that whp the induced partition of the graph achieves a modularity score near to $1$. Let $k$ be a large constant, let
$A_c=\{v_i : \Theta_i \in ((c-1)/k, c/k] \}$ and let $\cA=\{A_1, A_2, \ldots, A_k\}$ (so that we are radially partitioning the space into $k$ equal sectors).


\begin{center}
  \section{\hspace{-.3in} . \hspace{.1in}Modularity and partially observed graphs}
\label{sec.edge-sampling}
\end{center}
\vspace{-0.15in}
In this section, we first consider the percolated random subgraph $G_p$ of a graph $G$, where the edges are retained independently with probability $p$. We then briefly consider the more general case where the edge-probabilities may differ.  Finally, we discuss the `vertex-sampling' model, where we sample a random $k$-set of vertices and observe the corresponding induced subgraph.

\smallskip

\needspace{3\baselineskip}
\begin{center}
{\bf The percolated random graph model}
\end{center}
\vspace{-0.1in}
Given a graph $G$ and $0<p \leq 1$, let the {\em percolated} random graph $G_p$ be the random subgraph of~$G$ on the same vertex-set, obtained by considering each edge of $G$ independently, keeping it in the graph with probability~$p$, and otherwise deleting it. We assume that we can see the graph $G_p$ but not the underlying graph $G$, so that $G$ is `partially observed'. What can we infer about $G$? 
The binomial random graph $\Gnp$ (discussed in Section~\ref{sec.binom}) may be written as $(K_n)_p$, where $K_n$ is the $n$-vertex complete graph.  Let $e(H)$ be the number of edges in a graph $H$, so that the expected number of edges in $G_p$ is $e(G) p$.

\begin{figure}[h!]
    \centering
    \includegraphics[scale=0.6]{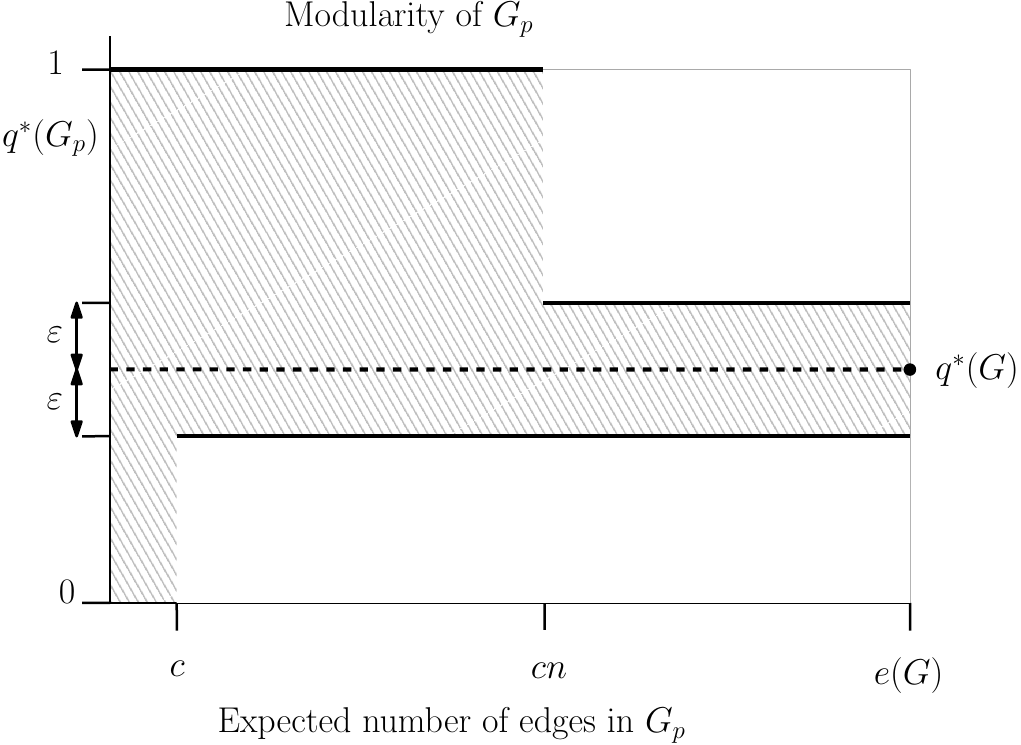}
    \caption{For any $\varepsilon>0$, there exists $c$ such that the modularity of the percolated graph $G_p$ lies in the region shown with probability at least $1-\varepsilon$. See Theorems~\ref{thm.obsmod} and~\ref{thm.moddiff}(a).}\label{fig.samplingobs}
\end{figure}

\begin{figure}[h!]
    \centering
    \vspace{5mm}
    \includegraphics[scale=0.6]{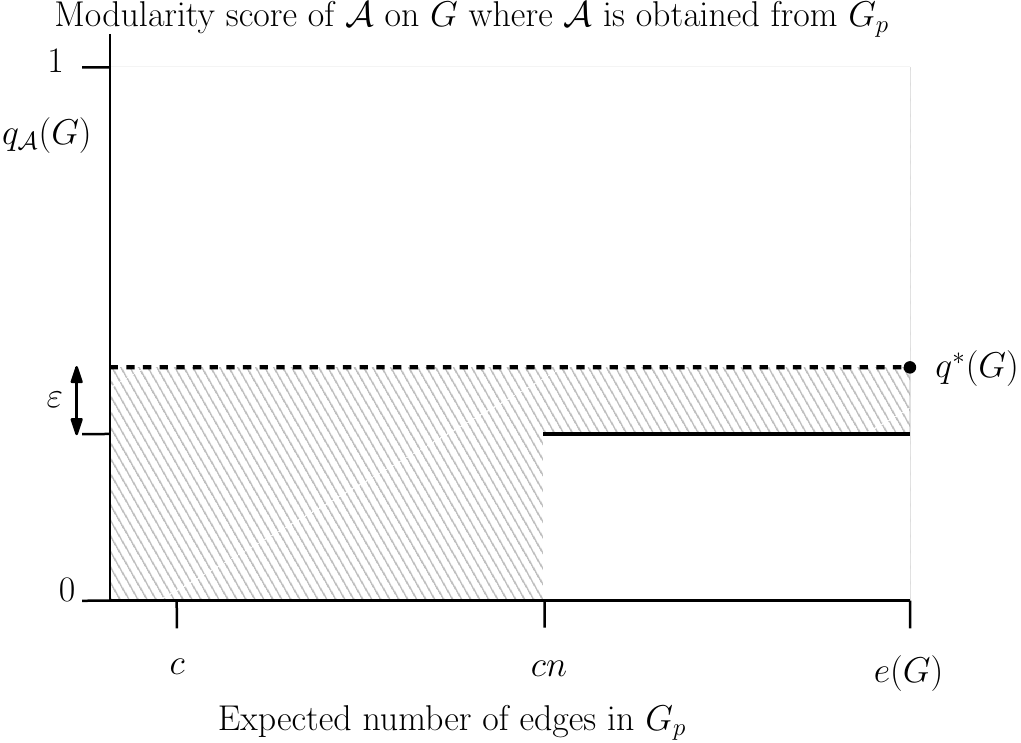}
    \caption{For any $\varepsilon>0$, there exists $c$ such that the partition $\cA$ constructed looking only at~$G_p$ has score $q_\cA(G)$ in the region shown, with probability at least~$1-\varepsilon$. See Theorem~\ref{thm.moddiff}(b).}\label{fig.samplingdiff}
\end{figure}

In this section we present two theorems from~\cite{sampling} that relate $\q(G_p)$ and $\q(G)$, as well the modularity score on $G$ of partitions constructed looking only at $G_p$; see  Fig.~\ref{fig.samplingobs} and Fig.~\ref{fig.samplingdiff}.
The first theorem concerns when we want to have $\q(G_p) > \q(G)-\eps$ with probability near~$1$, and it shows that this happens if the expected number of edges seen is large enough.
\begin{thm}\label{thm.obsmod}
Given $\eps>0$, there exists $c=c(\eps)$ such that the following holds. For each graph~$G$ and probability $p$ such that $e(G)p \geq c$, the random graph $G_p$ satisfies $\q(G_p) > \q(G)-\eps$ with probability $\geq 1-\eps$.
\end{thm}
\vspace{-0.1in}
\noindent
The proof of Theorem~\ref{thm.obsmod} in~\cite{sampling} shows that we may take $c(\eps)=\Theta(\eps^{-3}\log \eps^{-1})$, and an example graph involving a star and some disjoint edges shows that $c(\eps)$ must be 
$\Omega(\eps^{-1} \log \eps^{-1})$. 

The second theorem concerns when we want also to have $\q(G_p) < \q(G)+\eps$ with probability near 1, and it shows that this will happen if the expected average degree is sufficiently large.  
We see also that, in this case, finding a good partition $\cA$ for $G_p$ helps us to find a good partition $\cA'$ for $G$. 
Observe that the expected average degree in~$G_p$ is $2 e(G)p/v(G)$, where $v(G)$ is the number of vertices in $G$.

\needspace{4\baselineskip}
\begin{thm}\label{thm.moddiff}
For each $\eps>0$, there exists $c=c(\eps)$ such that the following holds.  Let the graph $G$ and probability $p$ satisfy $e(G)p/v(G) \geq c$.  Then, with probability $\geq 1-\eps$,\\
{\rm (}a{\rm )} the random graph $G_p$ satisfies $|\q(G_p) - \q(G)| < \eps$; and\\
{\rm (}b{\rm )} given any partition $\cA$ of the vertex-set, and seeing only~$G_p$, we can find in a linear number of operations a partition $\cA'$ with $q_{\cA'}(G) \geq q_{\cA}(G_p) - \eps$.
\end{thm}
\vspace{-0.1in}
\noindent
Here, part (b) says roughly that, given a good partition $\cA$ for $G_p$, we can quickly construct a good partition $\cA'$ for $G$. 
Using also part (a), we may see that, if the partition $\cA$ for $G_p$ satisfies $q_\cA(G_p) \geq \q(G_p) -\eps$ with probability at least $1-\eps$, then the partition $\cA'$ satisfies $q_{\cA'}(G) > \q(G) - 2\eps$ with probability at least $1-2\eps$. (See Section~3 of~\cite{sampling} for the greedy amalgamating algorithm used to construct $\cA'$ from $\cA$, where every part of $\cA'$ is a union of parts of $\cA$.)

The proof of Theorem~\ref{thm.moddiff} in~\cite{sampling} shows that we may take $c(\eps)=\Theta(\eps^{-3} \log \eps^{-1})$. An example with $G=K_n$ (so that $G_p$ is $\Gnp$ and $\q(G)=0$) shows that $c(\eps)$ must be at least $\Omega(\eps^{-2})$.  Taking $G=K_n$ again and applying part\,(a) of Theorem~\ref{thm.moddiff} shows that $\q(\Gnp) \to 0$ if $np \to \infty$, and so we have another proof of Theorem~\ref{thm.usER}\,(c).

The assumption in Theorem~\ref{thm.moddiff} that the expected average degree is large is much stronger than the assumption in Theorem~\ref{thm.obsmod}, but still $e(G_p)$ may whp be much smaller than $e(G)$. If we go further, and assume that
at most an $\eps$-proportion of edges are missed (that is, $|E \setminus E'| \leq \eps |E|$), then deterministically we have $|\q(G')-\q(G)| \leq 2 \eps$, by Theorem~\ref{lem.robustness}.

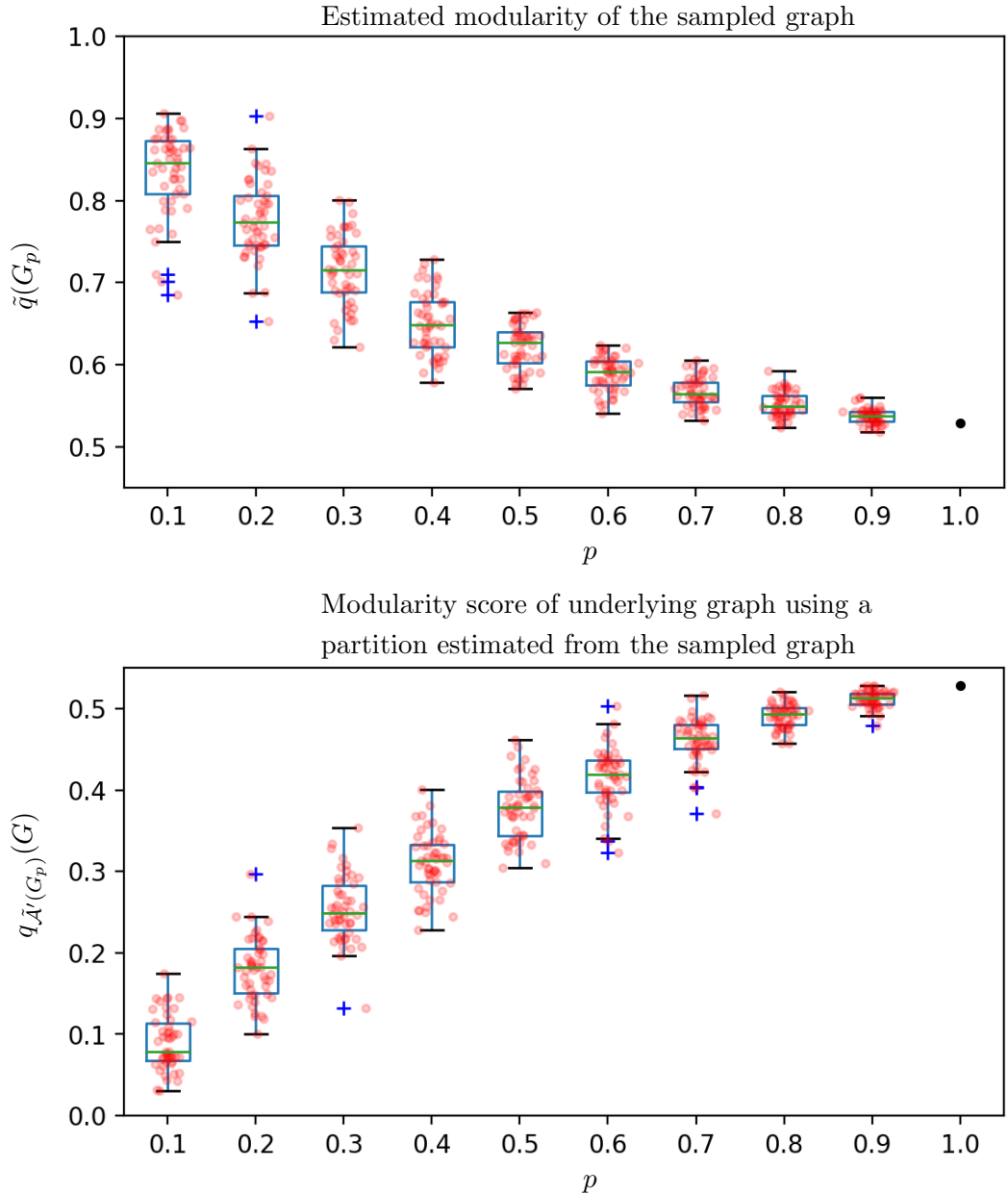
\begin{figure}
\begin{picture}(170,450) 
	\put(20,250){ 
		\includegraphics[width=0.95\textwidth]{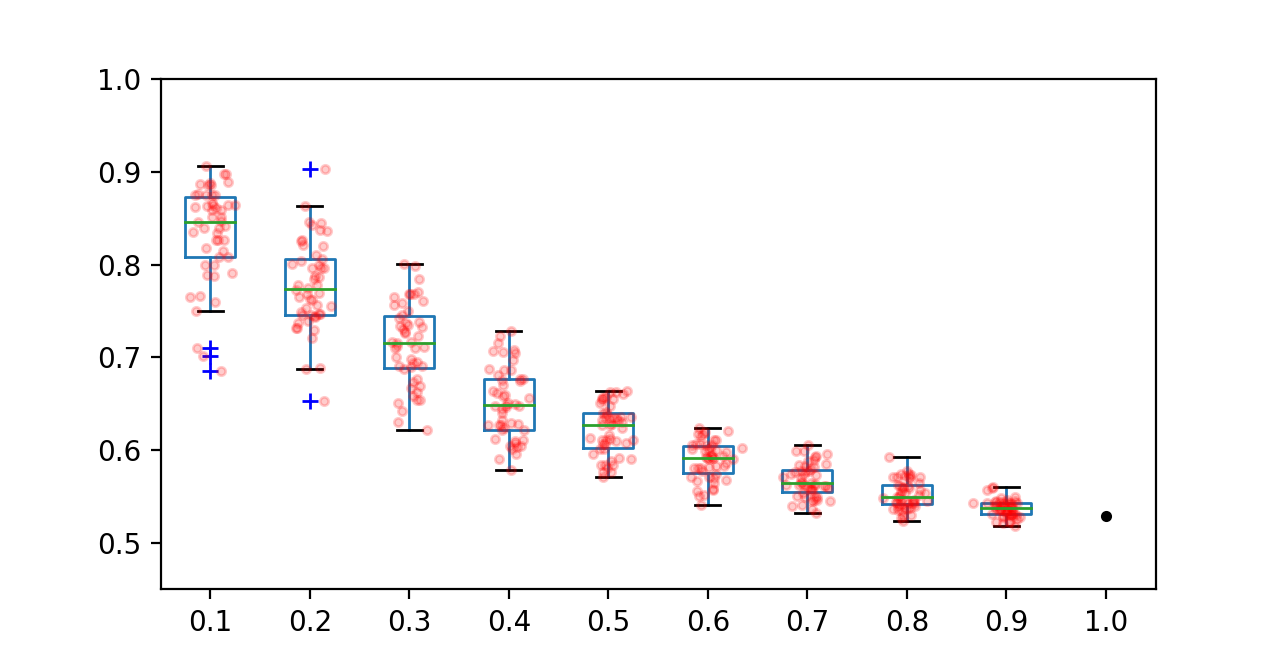}
		\put(-435,95){\rotatebox{90}{$\tilde{q}(G_p)$}} 
	    \put(-313,207){Estimated modularity of the sampled graph} 
	    \put(-210,-3){$p$}
		}
	\put(20,3){
	\includegraphics[width=0.95\textwidth]{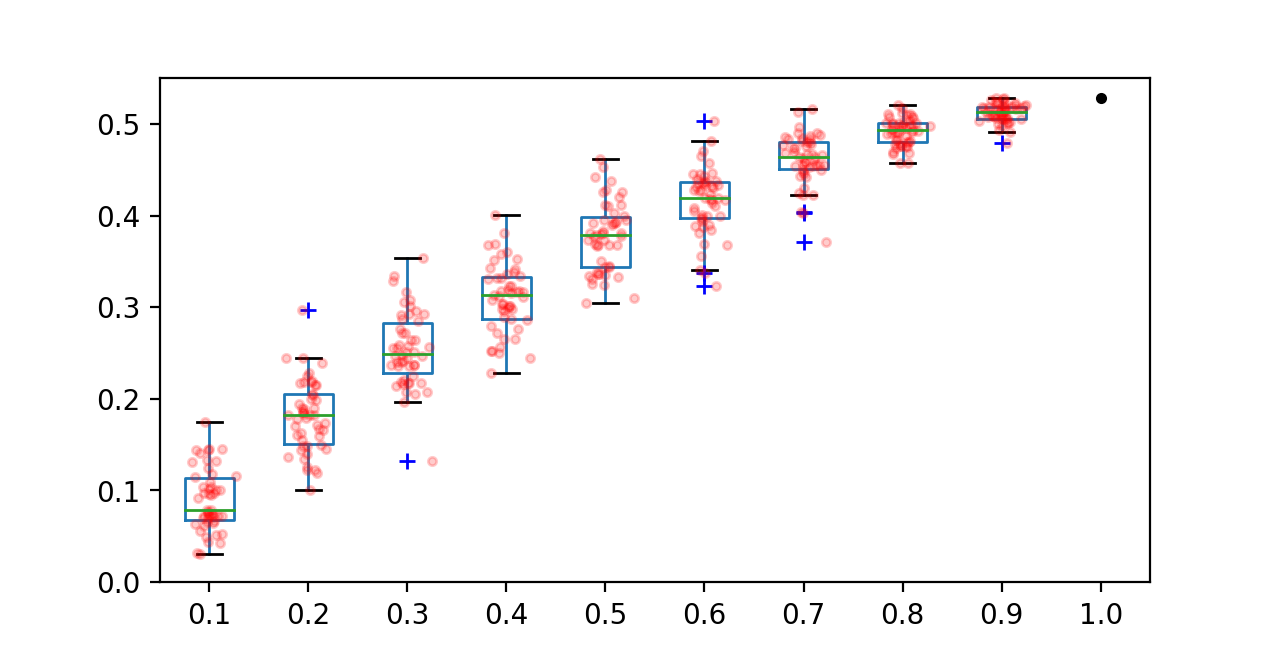}
		\put(-435,95){\rotatebox{90}{$q_{\tilde{\cA}'(G_p)}(G)$}} 
	    \put(-313,223){Modularity score of underlying graph using a}
	    \put(-313,207){partition estimated from the sampled graph}
	    \put(-210,-3){$p$}
		}
\end{picture}
\caption{Simulation results, based on the dolphin social network $G$, with $\q(G) \approx 0.529$.  In the upper part, each point corresponds to the estimated modularity $\tilde{q}(G_p)$ of an instance of the sampled graph~$G_p$.
(The noise in the $x$-axis is to allow one to see the points.)
In the lower part, for each sampled graph $G_p$, the corresponding point shows the score $q_{\tilde{\cA}'}(G)$, where $\tilde{\cA}$ is the highest scoring partition found for $G_p$, and~$\tilde{\cA}'$ is the partition for $G$ modified from $\tilde{\cA}$ as in Theorem~\ref{thm.moddiff}(b).
}
\label{fig.sims}
\end{figure}

Suppose that we take as the underlying graph $G$ the dolphin social network (see~\cite{lusseau2003emergent}), with 62 vertices and 159 edges. It is known that $\q(G) = 0.529$ to three decimal places (see~\cite{nphard}). In the upper part of Fig.~\ref{fig.sims}, each point corresponds to the estimated modularity $\tilde{q}(G_p)$ of an instance of the sampled graph $G_p$.
For each edge-probability $p=0.1, 0.2, \ldots, 0.9$, the graph $G_p$ was sampled 50 times. For each sampled graph $G_p$, the plot shows the maximum modularity score of the partitions output by 200 runs of both the Louvain~\cite{louvain} and Leiden~\cite{traag2019leiden} algorithms. 
The lower part of Fig.~\ref{fig.sims} examines, for each random instance of $G_p$, how well the modularity maximising partition found for~$G_p$ yields a good partition for the underlying graph $G$. In detail, for each sampled graph~$G_p$, the plot shows the score $q_{\tilde{\cA}'}(G)$, where $\tilde{\cA}$ is the highest scoring partition on $G_p$ found in the runs of the Louvain and Leiden algorithms, and $\tilde{\cA}'$ is the partition modified from $\tilde{\cA}$ as in Theorem~\ref{thm.moddiff}(b). For more details and simulations on a larger network, see~\cite{sampling}.

Let us discuss briefly how to prove Theorems~\ref{thm.obsmod} and~\ref{thm.moddiff}. A key preliminary result for both proofs is the `fattening lemma' below. 
Given a graph $G$ and a number $\eta$ with $0<\eta \leq 1$, we say that a vertex-partition $\cA$ is $\eta$-\emph{fat} if each part has volume at least $\eta \, \vol(G)$. 
There is a natural greedy algorithm (see~\cite{sampling}) which, given a graph $G$, a number $\eta$ with $0<\eta \leq 1$, and a vertex-partition $\cB$, amalgamates some parts in $\cB$ to yield an $\eta$-fat partition $\cA=\cA(G,\eta,\cB)$ with modularity score at least about that of $\cB$. Note that $|\cA| \leq |\cB|$. 
\begin{thm}[The fattening lemma] \label{lem.nosmall2}
For each graph $G$ and each number $\eta$ with $0<\eta \leq 1$, there is an $\eta$-fat partition $\cA$ of $V(G)$ such that $q_{\cA}(G) > \q(G) - 2\eta$. Indeed, given any partition $\cB$ of~$V(G)$, the greedy amalgamating algorithm uses a linear number of operations and constructs an $\eta$-fat partition $\cA =\cA(G,\eta,\cB)$ such that $q_{\cA}(G) > q_{\cB}(G) - 2\eta$.
\end{thm}
\vspace{-0.1in}
\noindent
The main step in the above algorithm is a number--bipartitioning problem: we have a set of parts $B_i$ from $\cB$ which we want to split into two sets, each with its sum of volumes at least $\eta\, \vol(G)$.
The constant~2 in Theorem~\ref{lem.nosmall2} is best possible (see~\cite{sampling}).
For comparison, recall Theorem~\ref{prop.1overk} 
concerning the maximum modularity score $\qk(G)$ for a partition of $G$ with at most~$k$ parts.
Observe that if $\cA$ is an $\eta$-fat partition for $G$ and $k = \lfloor 1/\eta \rfloor$, then $|\cA| \leq k$, and so $q_{\cA}(G) \leq \qk(G)$.  
\smallskip

To prove Theorem~\ref{thm.obsmod}, we note first that we may assume that
$\q(G) \geq \eps$, and let $\eta =\eps/4$.  By Theorem~\ref{lem.nosmall2}, there is an $\eta$-fat partition $\cA$ for $G$ with $q_\cA(G) \geq \q(G)-2\eta$. Since $|\cA| \leq \eta^{-1}$, we have $q^D_\cA(G) \geq \eta\,$, and also $q^E_\cA(G) = q_\cA(G) + q_\cA^D(G) \geq 3\eta$. 
Since $\cA$ is $\eta$-fat, and both $q^E_\cA(G)$ and $q^D_\cA(G)$ are bounded away from 0, we can use Chernoff probability bounds to show that $q^E_\cA(G_p)$ is unlikely to be much smaller than $q^E_\cA(G)$, and similarly that $q^D_\cA(G_p)$ is unlikely to be much larger than $q^D_\cA(G)$. The theorem follows.

\smallskip

The proof of Theorem~\ref{thm.moddiff} is more involved. The rough idea is that we can use the fattening lemma to bound the probability that a vertex-partition behaves badly by the probability that a fat vertex-partition behaves similarly badly. We can then use probabilistic methods (including Chernoff bounds, as above) to handle fat partitions.

\smallskip
\needspace{5\baselineskip}
\begin{center}
{\bf Percolation allowing distinct edge probabilities}
\end{center}
\vspace{-0.1in}
Up to now we have considered the model $G_p$, where each edge from an underlying graph $G$ appears independently with probability $p$, but we may also consider a model in which each edge $uv$ in $G$ appears in the observed graph with probability $p_{uv}$. Again, under certain conditions, the modularity of this graph converges whp to a deterministic quantity, and this is the modularity of a weighted graph.

Let $V$ be a non-empty vertex-set, and let $w: V\times V\to \R$ satisfy $w_{uv} = w_{vu} \geq 0$ for all vertices~$u$ and~$v$. We call $w$ a \emph{weight function} on $V^2$. 
Define the (weighted) degree of a vertex $u$ by setting $d^w_u =\sum_v w_{uv}$.
Similarly, define the (weighted) volume of a vertex-set $X$ by $\vol_w (X) = \sum_{u \in X} d^w_u$, and (corresponding to~$e(X)$) let $e_w(X)= \tfrac12 \sum_{u,v \in X} w_{uv}$, with $m_w=e_w(V)$ (corresponding to $m$).

Assume that $w$ is not identically zero, so that $m_w >0$. 
For a given partition $\cA$ of $V$, 
the \emph{modularity score} of $\cA$ on $w$ is 
\begin{equation}\label{eq.defn_w}
q_\cA(w) = 
\frac{1}{2m_w} \sum_{A\in \cA} \sum_{u,v \in A} \bigg(w_{uv} -\frac{d^w_u d^w_v}{2m_w}\bigg)\,, 
\end{equation}
and the \emph{modularity} of $w$ is $q^*(w)=\max_{\cA} q_{\cA}(w)$.
As in the unweighted case, $0 \leq \q(w)<1$ and we may ignore vertices with degree 0. Also, if $w$ is identically 0 we set $q_{\cA}(w)=0$ and $\q(w)=0$.  
If $w$ is $\{0,1\}$-valued and $w_{vv}=0$ for each vertex $v$, we recover the usual modularity score and modularity value  
for the corresponding graph.
The values $q_{\cA}(w)$ and $\q(w)$ are unchanged under a re-scaling of~$w$, and so we may assume that $0 \leq w_{uv} \leq 1$ for each $u,v \in V$, and in this case we call~$w$ a \emph{probability weight function}.

Given a probability weight function $w$, let $G_{w}$ be the random (unweighted) graph obtained by considering each edge independently, and including each edge $uv$ with probability $w_{uv}$.

\begin{thm}\label{thm.obsmod_wb} 
Given $\eps>0$, there exists $c=c(\eps)$ such that the following holds.
If the probability weight function $w$ on $V^2$ satisfies $e_w(V) \geq c$, then with probability at least $1-\eps$, the random graph~$G_{w}$ satisfies $\q(G_{w}) > \q(w) - \eps$. 
\end{thm}
\vspace{-0.1in}
\noindent

\needspace{3\baselineskip}
\begin{thm}\label{thm.moddiff_wb}
Given $\eps>0$, there exists $c=c(\eps)$ such that the following holds.
If the probability weight function $w$ on $V^2$ satisfies $e_w(V) \geq c\, |V|$, then with probability at least $1-\eps$,\\ 
{\rm (}a{\rm )} the random graph $G_{w}$ satisfies $|\q(G_{w}) - \q(w)| <\eps$; and\\
{\rm (}b{\rm )}  given any partition $\cA$ of the vertex-set, in a linear number of operations (seeing only~$G_w$) the greedy amalgamating algorithm finds a partition $\cA'$ with 
\mbox{$q_{\cA'}(w) \geq q_{\cA}(G_{w}) \!-\!\eps$}.
\end{thm}
\vspace{-0.1in}
\noindent
In the next section, Theorem~\ref{thm.moddiff_wb} is used to establish the modularity score of the stochastic block model.

\smallskip

\begin{center}
{\bf Estimating modularity by vertex sampling (`parameter estimation')}
\end{center}
\vspace{-0.1in}
Another common model of partially observing an underlying graph $G$ is to sample a vertex-subset~$U$ of constant size $k$, and observe the corresponding induced subgraph~$G[U]$.
A graph parameter $f$ is \emph{estimable} (or \emph{testable}) for a class $\mathcal G$ of graphs
(see~\cite{BCLSVjournal} and \cite{lovasz2012large}) if, for every $\eps>0$, there exists $k = k(\eps)$ such that, if $G$ is a graph in $\mathcal G$ with at least $k$ vertices, then $\pr(|f(G)-f(G[X])| >\eps)<\eps$ when $X$ is a uniformly random $k$-subset of $V(G)$. The following result (\cite[Theorem~1.5]{sampling}) shows
that the modularity value~$\q(G)$ is estimable for dense graphs~$G$, but not more generally.
\begin{thm}\label{thm.modest} 
{\rm (}a{\rm )} For any fixed $\rho$ with $0<\rho<1$, modularity is estimable for graphs with density at least~$\rho$.\\
{\rm (}b{\rm )} For any given function $\rho(n) = o(1)$, modularity is not estimable for $n$-vertex graphs $G$ with density at least $\rho(n)$ (even if we assume that $G$ is connected).
\end{thm}
\vspace{-0.1in}
\noindent
To understand the proof of Theorem~\ref{thm.modest}  we need some definitions.
If $G$ is a graph and $S, T \subseteq V(G)$, then $e_G(S,T)$ is the number of ordered pairs $(s,t) \in S \times T$ for which there is an edge in $G$ between~$s$ and $t$; note that
$S$ and $T$ may overlap. If $G$ and $G'$ are graphs with the same vertex-set $V$, then the \emph{cut distance} is 
\[ d_\square(G,G') = |V|^{-2} \max_{S,T}|e_G(S,T)- e_{G'}(S,T)|\,.\]
Given a graph $G$ and a positive integer $b$, we let $G(b)$ denote the $b$-\emph{blow-up} of $G$, where each vertex of $G$ is replaced by $b$ independent copies of itself; thus, $v(G(b))=b\,v(G)$ and $e(G(b))= b^2\, e(G)$.

\needspace{4\baselineskip}
By \cite[Theorem~15.1]{lovasz2012large}, $\q(G)$ is estimable for $\mathcal G$ if and only if the following two conditions hold:\\
(i) If graphs $G_n$, $G'_n$ in $\mathcal G$ have the same vertex-set and $d_\square(G_n, G'_n)\rightarrow 0$ as $n \to \infty$, then $\q(G_n)-\q(G'_n)\rightarrow~\!0$ as $n \to \infty$.\\
(ii) For every graph $G$ in $\mathcal G$, $\q(G(b))$ has a limit as $b \rightarrow \infty$.

\vspace{-0.1in}
\noindent
It is straightforward to show that (ii) holds, and in fact $\q(G(b)) \to$ 
$\sup_b \q(G(b))$ as $b \rightarrow \infty$.  By considering separately the edge contribution and degree tax, much as in the proof of Theorem~\ref{lem.robustness}, we can show that (i) holds for graphs with density at least $\rho$.

Finally, let us see that (i) fails (badly) in the non-dense case, even for connected graphs.  Let $0 \leq \rho(n)<1$, and let $\rho(n) \to 0$ arbitrarily slowly as $n \to \infty$;  we may assume that $\rho(n) n/\log n \to~\infty$. 
Then there are connected graphs $G_n$, $G'_n$ on the vertex-set $[n]$  with density at least $\rho(n)$, such that each has density $o(1)$, so $d_\square(G_n, G'_n)\rightarrow 0$, but $\q(G_n) \to 1$ and $\q(G'_n) \to 0$.
We may take $G_n$ to be a collection of disjoint cliques of size about $2 \rho(n) n$, joined by edges to form a path.  For $G'_n$, we may consider a binomial random graph $G_{n,\rho}$\,, which whp is connected, has density between $\rho/2$ and $2 \rho$, and by Theorem~\ref{thm.usER}(c) has modularity~$o(1)$.


\smallskip

\begin{center}
  \section{\hspace{-.3in} . \hspace{.1in}
  Modularity and stochastic block models} 
  \label{sec.stoch-block}
\end{center}
\vspace{-0.15in}
We focus mostly on balanced stochastic block models.
Let  $n, k \geq 2$, and let $0 \leq q \leq p \leq 1$: we define the corresponding stochastic block model $G_{n, k,p,q}$ as follows.  The vertex-set $V$ is $[n]$, and the random variables $\sigma_v$ for $v \in V$ are independent, with each uniformly distributed on $[k]$. This yields $k$ random blocks or communities $V_i=\{v \in V: \sigma_v=i\}$, each with expected size $n/k$, which form the `planted partition' of $V$. Edges appear independently, with probability~$p$ within blocks and with probability~$q$ between blocks. Typically we are not told the planted partition.

For each fixed $k \geq 2$, as long as the average degree tends to $\infty$, the planted partition whp yields a modularity score which is close to optimal.

\begin{thm} \label{thm.SBMk}
Let $k \geq 2$, and let $p=p(n)$ and $q=q(n)$ satisfy $0  \leq  q \leq p \leq 1$ and $np \to \infty$ as $n \to \infty$.
Then   
\[ \q(G_{n, k, p, q}) = 
\frac{(p-q) \, (1 - 1/k)}{ p + (k-1)q} +o(1) \;\;\; \mbox{whp},\]
and whp the planted partition has this modularity value.
\end{thm}
\vspace{-0.1in}
\noindent
A direct proof of this result is given in~\cite{sampling}, using ideas from~\cite{koshelev2023modularity}. Given the planted partition, we construct a weighted graph with edge-weights $a$ within communities and $b$ between communities; as we saw in Section~\ref{sec.edge-sampling}, there is a natural way to define the modularity of such a weighted network.
Using Theorem~\ref{thm.moddiff_wb}, we can show that the stochastic block model has modularity value very close to that of this underlying weighted graph. Thus it remains only to determine the modularity value of a simple deterministic weighted graph, and this is quite straightforward. 

Let us now consider a different way to prove a result like Theorem~\ref{thm.SBMk}, at least in the case when the ratio $p/q$ is fixed.
We saw above that the planted partition whp yields a close-to-optimal modularity value, but is it close to an optimal partition?
The following result, Theorem~\ref{thm.SBMrecovery}  
(see~\cite{bickel2015correction} and the references therein), answers this question. We can thus give another proof of a result like Theorem ~\ref{thm.SBMk}, by using Theorem~\ref{thm.SBMrecovery} and the robustness result Theorem~\ref{lem.robustness} and calculating the likely modularity score of the planted partition.

\begin{thm}\label{thm.SBMrecovery}
Let $a>b>0$, and let $\omega=\omega(n) \to \infty$ as~$n \to \infty$. Let $\cA^*$ be a modularity-optimal partition of $G_{n,{k},p,q}$.\\
If $p=a\, \omega/n$ and $q=b\, \omega /n$, then whp $\cA^*$ and the planted partition differ on $o(n)$ vertices.\\
If $p=a\, (\omega \log n)/n$ and $q=b\, (\omega \log n) /n$, then whp $\cA^*$ is exactly the planted partition. 
\end{thm}
\noindent We next define a generalised $k$-block stochastic block model. Let $k \geq 2$, let $\pi=(\pi_1,\pi_2, \ldots, \pi_k)$ with~$\pi_i \geq 0$ and $\sum_i \pi_i=1$, let~$\mathbf{P}=[p_{ij}]$ be a $k\times k$ matrix with $0 \leq p_{ij} \leq 1$, let $n\geq 2$, and let $\rho=\rho(n)$ with $0<\rho \leq 1$.  We define $G_{n,\pi, \mathbf{P}, \rho}$ as follows. 
The vertex-set~$V$ is~$[n]$, and for~$v\in V$, we let the random variables $\sigma_v$ be independent and identically distributed, taking values in $[k]$ with probability $\pi_i$ of taking value $i$. This yields the $k$-block `planted partition', where $V_i =  \{ v \in V :  \sigma_v = i\}$ with expected size $\pi_i n$. The edges within each block $V_i$ appear with probability $\rho\,  p_{ii}$ (following~\cite{zhao2012consistency}, we allow loops), and the edges between blocks $V_i$ and $V_{j}$ appear with probability $\rho\, p_{ij}$.

When the probabilities within blocks and between blocks can depend on the labels of the blocks, the problem of recovering communities is subtle.
A balanced $3$-block example is given in~\cite{bickel2009nonparametric}, in which a partition that puts two blocks together has a higher modularity score than the planted partition. In this example, $p_{11}=0.06$, $p_{22}=0.12$, $p_{33}=0.66$, $p_{12}=p_{23}=0.04$, $p_{13}=0$, and~$\rho=1$.
Let $\cA$ be the planted partition (with three parts), and let $\cA'$ be the bipartition which merges blocks~$1$ and~$2$. Then whp $q_\cA(G) \sim 0.84-0.54=0.3$ and $q_{\cA'}(G)\sim 0.92-0.58= 0.34$, where we recall that $\sim$ indicates that equality holds with a $(1+o(1))$ factor, and where the two terms are the edge contribution and degree tax. 
In this example, the minimum within-block
probability is higher than the maximum between-block 
probability.
For the following theorem, see~\cite{bickel2009nonparametric}, \cite{bickel2015correction} and~\cite{zhao2012consistency}.
\begin{thm}\label{thm.SBMrecovery_general}
Let $G_{n,\pi, \mathbf{P}, \rho}$ be a generalised $k$-block stochastic block model as above, and let~$\cA^*$ be a modularity-optimal partition of $G_{n, \pi,  \mathbf{P}, \rho}$. Let $\bar{p}=\sum_{i,j} \pi_i\pi_j p_{ij}$, and suppose that, for all $i\neq j$,
\vspace{-0.02in}
\begin{equation}\label{eq.zhao_orig} p_{ii} > \big(\sum_{h} \pi_h\, p_{ih}\big)^2/\bar{p}  \;\; \mbox{ and } \;\; p_{ij} < \big(\sum_{h} \pi_h p_{ih}\,\big)\big(\sum_h \pi_h\, p_{jh}\big)/\bar{p}\,.
\vspace{-0.02in}
\end{equation}
Let $\omega=\omega(n) \to \infty$ as $n \to \infty$.
If $\rho =\omega/n$  then whp $\cA^*$ and the planted partition differ on $o(n)$ vertices, and if $\rho = (\omega \log n) /n$ then whp $\cA^*$ is exactly the planted partition. 
\end{thm}
\noindent
Condition~\eqref{eq.zhao_orig} on $\pi$ and $\mathbf{P}$ is equivalent to the following inequalities for the random graph~$G$ constructed from $G_{n, \pi, \mathbf{P}, \rho}$ by removing loops:
\begin{equation}\label{eq.zhao_condition_exp}
\frac{\E[e(V_i)]}{\E[e(G)]} > \frac{\E[\vol(V_i)]^2}{\E[\vol(G)]^2} 
\;\;\; \mbox{ and } \;\;\; \frac{\E[e(V_i, V_j)]}{2\,\E[e(G)]} < \frac{\E[\vol(V_i)]\,\E[\vol(V_j)]}{\E[\vol(G)]^2}.\end{equation}
(Here $e(V_i)$ means $e_G(V_i)$, etc.) 

Returning to the balanced model $G_{n,k,p,q}$, we note that~\cite{bickel2015correction} and~\cite{zhao2012consistency} show that Theorem~\ref{thm.SBMrecovery} holds also for modified versions of the modularity function, including the Erd\H{o}s--R\'enyi modularity, where the~$d_ud_v/ 2m$ degree tax term is replaced by $2m/n$. The same papers extend Theorem~\ref{thm.SBMrecovery} to a `degree-corrected' version of the stochastic block model, using a quality function similar to modularity. %

{The lower bound on the modularity $\q(G_{n,p})$ in~Theorem~\ref{thm.growthRate} covers a wide range of probabilities~$p$, and has a stand-alone algorithmic proof. Recall that when $1/n \leq p \leq 1- c_0/n$, the algorithm \emph{Swap} whp finds a balanced bipartition achieving a modularity score of at least $\alpha\sqrt{(1\!-\!p)/np}$ (and the constant~$\alpha$ may be taken to be $\frac15$ in part of that range). For the special case when $p=c/n$, attractive results on contiguity between binomial random graphs and stochastic block models (see~\cite{banks2016information} and~\cite{distinguish}) yield another proof giving a better constant. 
\begin{thm} \label{thm.lowerSqrt2}
For each constant $c > 1$, 
$\q(\Gncntext) >  0.668 / \sqrt{c}$ \mbox{ whp}. 
\end{thm}
\vspace{-0.1in}
\noindent
This result may be compared to the value $\q(\Gncntext)\sim 0.97/\sqrt{c}$\,  predicted using spin-glass models (see~\cite{trulymodular}). We recall also the bound in~\cite{rybarczyksulkowska2025Gnp} that, for $c$ sufficiently large, whp $\q(\Gncntext)\geq 0.7632/\sqrt{c}$, which is stronger for large $c$ -- see~\eqref{eqn.b} for a corresponding lower bound.
To prove Theorem~\ref{thm.lowerSqrt2}, we can show that, for each $k \geq 2$, whp there is a balanced $k$-part partition with a modularity score of at least about $f(k)/\sqrt{c}$, for an explicit function~$f(k)>0$.
In detail, let $f(2)=\frac12$, and for $k \geq 3$ let $f(k) = \sqrt{2(k\!-\!1)\log(k\!-\!1)} /k$\, (recall that $\log$ denotes natural logarithm).
Then, given a constant $c>1$, for each $k \geq 2$ whp there is a balanced $k$-part partition $\cA_k$ such that
\begin{equation} \label{eq.fk}
q_{\cA_k}(\Gncntext) 
\geq f(k) / \sqrt{c} \: + o(1),
\end{equation}
and choosing $k=6$ parts here yields the constant given in Theorem~\ref{thm.lowerSqrt2}.

We next sketch the proof of~\eqref{eq.fk} when $k=2$, and then briefly discuss the case $k \geq 3$. Let $0<b \leq a \leq n$, let $p=a/n$ and $q=b/n$, and consider the stochastic block model $G_{n,k,p,q}$.\\

\vspace{-0.05in}

\needspace{2\baselineskip}
\noindent
{\centering\textbf{Balanced bipartitions, $k=2$}\,

}
\vspace{0.01in}

\noindent When $a$ and $b$ are close together, the planted bisection model $G_{n,2,p,q}$ is \emph{contiguous} with $\Gncntext$ where $c=(a+b)/2$ -- that is, events $A_n$ hold whp in $G_{n,2,p,q}$ if and only if they hold whp in $\Gncntext$. A result of Mossel, Neeman and Sly~\cite{distinguish} says precisely that the models are contiguous if and only if $(a-b)^2 \leq 2(a+b)$. 
It follows that, if we fix $c >1$ and let $a=c+\sqrt{c}$ and $b=c-\sqrt{c}$, 
then the models $\Gncntext$ and $G_{n,2,p,q}$ are (just) contiguous. It is thus sufficient to show that whp we have the desired bipartition in $G_{n,2,p,q}$. %

As usual, let $\omega=\omega(n) \to \infty$ as $n \to \infty$, with $\omega=o(n)$.  Also, let $V_1=\{v\!: \sigma_v= 1\}$ and $V_2=\{v\!: \sigma_v=2 \}$, and let $\cA$ be the partition into $V_1$ and $V_2$, the planted partition. We can use Chebyshev's inequality repeatedly to show that, for $i=1$ and $2$, whp we have $|V_i| = \frac12 n +  o(\sqrt{\omega n})$, 
$e(V_i) = \frac18 a n +  o(\sqrt{\omega n})$, and
$\vol(V_i) =  \frac14 (a+b) n +  o(\sqrt{\omega n})$. 
Since $e(G_{n,2,p,q}) = \frac14 (a+b)n + o(\sqrt{\omega n})$ whp, we have
\begin{equation}
q_{\cA}(G_{n,2,p,q}) = \frac{a}{a+b} - \tfrac12 + o(\sqrt{\omega/n}) = \frac{a-b}{2(a+b)} + o(\sqrt{\omega/n}) \;\; \mbox{ whp},
\end{equation}
and so, by our choice of $a$ and $b$, whp $q_{\cA}(G_{n,2,p,q}) = 1/(2 \sqrt{c}) + o(\sqrt{\omega/n})$.
Further, whp the number of isolated vertices in $V_i$ is $\frac12 n e^{-c}+  o(\sqrt{\omega n})$, for $i=1,2$.
We may shuffle isolated vertices in a partition without changing the modularity, so whp we may modify $\cA$ to give a balanced partition, as required.\\


\needspace{3\baselineskip}
\noindent
{\centering\textbf{Balanced $k$-part partitions for $k\geq 3$.}

}
\vspace{0.05in}

\noindent By~\cite[Theorem~1]{banks2016information}, for $c=(a+(k-1)b)/k$, the models $\Gncntext$ and $G_{n,k,p,q}$ are contiguous if $(a - b)^2< 2ck^2\, \ln (k\!-\!1) \,/(k\!-\!1)$.  
Let $a=c+x\sqrt{c}$ and $b =c - x\sqrt{c} / (k-1)$, where $0<x< f(k)$, and where we take $x$ near the upper bound.
We now argue much as when $k=2$.  We see that~$\Gncntext$ and~$G_{n,k,p,q}$  
are contiguous, and use Chebyshev's inequality repeatedly to check that, for the planted partition $\cA$, 
$q_{\cA}(G_{n,k,p,q}) =  f(k)/\sqrt{c} +o(1) $ whp; 
finally, we shuffle isolated vertices to obtain a balanced partition.

\smallskip

\begin{center}
    \section{\hspace{-.3in} . 
    \hspace{.03in}
    The modularity of random graphs embeddable in a surface} \label{sec.furthermodels}
\end{center}
\vspace{-0.15in}
The modularity of a random planar graph is close to 1 whp.  
More generally, let us consider the modularity of a random graph embeddable in a given surface $S$.
Let $\cE^S$ be the class of graphs embeddable in $S$. We then have the following theorem.

\begin{thm} \label{prop.fixedS}
Let $S$ be any surface, and let $R_n \in_u \cE^S$.  Then $\q(R_n) = 1-o(1)$ whp.
\end{thm}
\vspace{-0.1in}
\noindent
How can we prove this?
The \emph{tree-width} ${\rm tw}(G)$ of a graph $G$ is a measure of how similar it is to a tree (see for example~\cite{kloks1994treewidth}); recall that $\Delta(G)$ is the maximum vertex-degree. 
Consider the following deterministic result
\cite[Corollary~1.12]{treelike}.
For $m = 1, 2, \ldots,$ let $G_m$ be a graph with $m$ edges. If ${\rm tw}(G_m) \cdot \Delta(G_m) = o(m)$, then
$\q(G_m) \to 1$ as $m \to \infty$.
But ${\rm tw}(R_n) = O(\sqrt{n})$, whp $\Delta(R_n) = O(\log n)$, and whp $e(R_n) = \Theta(n)$, and so  Theorem~\ref{prop.fixedS} follows.
Theorem~\ref{prop.fixedS} can also be deduced easily from the following result of Laso\'n and Sulkowska~\cite{lasonsulkowska2023modularity}. 
\begin{thm}
    \label{lem.minorfree} 
 Let the class $\cA$ of graphs have an excluded minor, and suppose that the maximum degree of any graph $G_m$ in $\cA$ with $m$ edges is $o(m)$. Then whp $\q(G_m)=1-o(1)$. 
\end{thm}
\vspace{-0.1in}
\noindent
To deduce Theorem~\ref{prop.fixedS} from Theorem~\ref{lem.minorfree}, let $\eps(m) = 1/\log m$, say, and let $\cB$ be the class of graphs $G \in \cE^S$ for which $\Delta(G) \leq \eps(m)\, m$, where $m=e(G)$.  Then  Theorem~\ref{lem.minorfree} applied to $\cB$, together with the results noted above that whp $\Delta(R_n) = O(\log n)$ and $e(R_n) = \Theta(n)$, again show that $\q(R_n) = 1-o(1)$ whp.

We can extend Theorem~\ref{prop.fixedS} by considering order-dependent surfaces (see~\cite{mcdiarmid2023random}). Let $g=g(n) \geq 0$, and let $\cE^g$ be the class of graphs $G$ such that, if $G$ has order $n$, then $G \in \cE^S$ for some surface $S$ with Euler genus at most $g(n)$. 
\begin{thm} \label{prop.Eg}$\mbox{ }$\\
{\rm (i)} Let $g(n) = o(n/\log^3n)$ and let $R_n \in_u \cE^g$; then $\q(R_n) = 1-o(1)$ whp.\\
\noindent {\rm (ii)}  Let $g(n)=o(n/\log n)$ be non-decreasing, and let $R_n \in_u \cE^g$; then there is a subsequence $1 \leq n_1<n_2< \cdots$ for which $\q(R_{n_k})=1-o(1)$ whp.
\end{thm}
\vspace{-0.1in}
\noindent
To prove part (i), we can show that $\Delta(R_n)$ is still $O(\log n)$ whp (see~\cite[Theorem 8(a)]{mcdiarmid2023random}).
Call a set $F$ of edges in a graph $G$ a \emph{planarising edge-set} if $G \setminus F$ is planar.
A result of Djidjev and Venkatesan (see~\cite[Lemma 27]{mcdiarmid2023random}) shows that, for any $n \geq 2$ and $h \geq 0$, any connected graph $G \in \cE^h_n$ 
has a planarising edge-set of size at most $4\sqrt{h(n+h) \Delta(G)}$.
So for $R_n$, whp we have a planarising edge-set $F_n$ of size $o(n)$, and so $|F_n|/e(R_n) =o(1)$.  For the planar graph $R_n^- = R_n \setminus F_n$, whp $\q(R_n^-) = 1+o(1)$, by the proof of Theorem~\ref{prop.fixedS}, and by the robustness result Theorem~\ref{lem.robustness}, we have $\q(R_n) = \q(R_n^-) +o(1) = 1 + o(1)$ whp.

We can prove part (ii) similarly; for by~\cite[Theorem 8(b)]{mcdiarmid2023random}, there are a subsequence $n_1<n_2< \cdots$ and a constant $c$ such that $\Delta(R_{n_k}) \leq c \log(n_k)$ whp (as $k \to \infty$).

\begin{center}
    \section{\hspace{-.3in} . \hspace{.1in} Some known modularity values}\label{sec.modZoo}
\end{center}
\vspace{-0.25in}
We conclude with a list of some known modularity values of graphs, sorted into four classes, those with modularity near to 1, 
modularity bounded away from 1 and from 0, modularity near to 0, and modularity exactly 0.
Knowing the modularity for classes of graphs may help us to understand the behaviour of the modularity function. We list results on the modularity of the binomial random graph $\Gnp$\,: 
as we saw earlier, there are similar results for the Erd\H{o}s--R\'enyi random graph $G_{n,m}$. 
We always use $n$ for the number of vertices and $m$ for the number of edges.
In the second column we give extra conditions, if these are necessary, and in the last column we refer to a result in this chapter or to the appropriate reference. \\

\footnotesize
$
\renewcommand{\arraystretch}{1.25}
\begin{array}{|l | l | l |  l | l}
\cline{1-1}
\mbox{\textbf{maximally modular}}\\ 
    \cline{1-4} 
\mbox{cycle} &
    & \q(C_n)=1-2/\sqrt{n}\,(1+o(1))
    &\mbox{\cite{nphard}}\\ 
    &&&&\\
\mbox{tree}  	
    & \Delta(T_n)=o(n) & \displaystyle \q(T_n) \geq 1- 2(2\Delta/n)^{1/2} 
    & \mbox{\cite{treelike}} \\
    &&&&\\
\mbox{tree-like (low tree-width)}
    & \Delta(G_m) {\rm tw}(G_m)=o(m) 
    & \displaystyle \q(G_m) \geq 1- 2(({\rm tw} +1)\Delta/m)^{1/2} 
    &\mbox{\cite{treelike}}  \\
    &&&&\\
\mbox{minor-free graphs}& 
\mbox{fixed forbidden minor }
    & \displaystyle \q(G_m) = 1 + o(1)
    &\mbox{Thm~\ref{lem.minorfree}}  \\
    &\& \; \Delta(G_m)=o(m)
    &&&\\
    &&&&\\
\mbox{binomial / Erd\H{o}s--R\'enyi}
    &n^2p\rightarrow \infty\; \& \;np\leq 1+o(1) 
    &\mbox{whp }\q(\Gnp)=1+o(1) 
    &\mbox{Thm~\ref{thm.usER}} \\
    &&&\\
\mbox{random hyperbolic}& 
    &\mbox{whp } \q(G_n) = 1 + o(1) 
    &\mbox{Thm~\ref{thm.hyperbolic}}\\
    &&&\\
\mbox{random 2-regular} 
    &%
    &\mbox{whp } \q(G_{n,2})=1-2/\sqrt{n}+{o(1)}
    &\mbox{Eqn~\eqref{eq.2_regular}}\\  
    &&&\\
\mbox{random geometric} 
    & 
    &\mbox{whp }  \q(G_n)=1 + o(1)
    &\mbox{\cite{davis2018consistency}}\\
    &&&&\\
\mbox{random planar}
    & 
    &\mbox{whp } \q(G_n) = 1 + o(1) 
    &\mbox{Thm~\ref{prop.fixedS}}\\
\mbox{(or any fixed surface)} 
    &&&&\\
\cline{1-4} \multicolumn{2}{l}{$\mbox{ }$}\\
%
%
    \cline{1-1}
\mbox{\textbf{critically modular}}\\
    \cline{1-4} 
\mbox{binomial / Erd\H{o}s--R\'enyi} 
    &\mbox{$np=c >1 \; \Rightarrow \; \exists\,\eps>0 \; :$} 
    & \mbox{whp } \eps <\q(\Gnp) < 1-\eps
    &\mbox{Thm~\ref{thm.usER}} \\
    &\mbox{$np=c \geq 1 \; $ } 
    &\mbox{whp } 0.668/\sqrt{c} <\q(\Gnp) < 3.06/\sqrt{c} 
    &\mbox{Eqn~\eqref{eqn.b}},\\
    &&&\mbox{Thm~\ref{thm.lowerSqrt2}}  \\
    &&&&\\
\mbox{random cubic} 
    &
    &\mbox{whp } 0.667 < \q(G_{n,3}) < 0.790   
    &\mbox{Thm~\ref{thm.cubic}}.   \\ 
    &&&&\\
\mbox{random $r$-regular} 
    &\mbox{for fixed }r\geq r_0  
    & \mbox{whp } 0.76/\sqrt{r} < \q(G_{n,r})<2/\sqrt{r} 
    & \mbox{\mbox{Thm~\ref{thm.rreg_larger} }} 
    \\
    &&&&\\
\mbox{preferential attachment} 
    & \mbox{$h\geq 2$ edges per step}
    & \mbox{whp } \Omega(1/\sqrt{h})=\q(G_{n}^h) = O(\sqrt{\log h}/{\sqrt{h}})  
    & \mbox{Thm~\ref{thm.PA}.} \\ 
    &&&&\\
\mbox{stochastic block model }
    & 0 \leq q\leq p \leq 1 \; \& \; np\rightarrow \infty 
    & \mbox{whp } \displaystyle \q({G_{n}}) =\frac{(p-q)(1-1/k)}{p+(k-1)q}+o(1) 
    &\mbox{Thm~\ref{thm.SBMk}}\\
    \mbox{($k$ balanced blocks)} &&&&\\\cline{1-4} \multicolumn{2}{l}{$\mbox{ }$}\\
\cline{1-1} 
\mbox{\textbf{minimally modular}}\\ 
    \cline{1-4}
\mbox{binomial / Erd\H{o}s--R\'enyi}
    & np\rightarrow \infty & \mbox{whp } \q(\Gnp)=o(1) 
    &\mbox{Thm~\ref{thm.usER}} \\ 
    \cline{1-4} \multicolumn{2}{l}{$\mbox{ }$}\\
    \cline{1-1} 
\mbox{\textbf{non-modular}}\\ \cline{1-4}
\mbox{complete}
    & 
    & \q(K_n)=0 
    & \mbox{\cite{nphard}}\\ 
    &&&&\\
\mbox{complete multipartite}
    & 
    & \q(K_{n_1,n_2,\ldots,n_k})=0 
    &\mbox{\cite{bolla2015spectral, vdense,majstorovic2014note}}\\
    &&&&\\
\mbox{nearly complete}
    & e(G)\geq\binom{n}{2}-n/2 
    & \q(G)=0
    &\mbox{\cite{vdense}}\\ 
    &&&&\\
\mbox{binomial / Erd\H{o}s--R\'enyi}
    & p\geq 1-c/n, \;\; c<1 
    & \mbox{whp } \q(\Gnp)=0 
    &\mbox{Thm~\ref{thm.Gnp}} \\  
    \cline{1-4}			
\end{array}
$

\normalsize
\needspace{10\baselineskip}


\footnotesize
\begin{thebibliography}{9}



\bibitem{mod2023universal}
  V.\ Agdur,  N.\ Kam{\v{c}}ev and F.\ Skerman,
  {Universal lower bound for community structure of sparse graphs},
  {arXiv:2307.07271}, (2023).





\bibitem{banks2016information}
	{J.~Banks, C.~Moore, J. Neeman and P.~Netrapalli}, {Information-theoretic thresholds for community detection in sparse networks},
	\emph{Conference on Learning Theory}
	(2016), {383--416}.


\bibitem{bickel2009nonparametric}
  {P.~J.\ Bickel  and A.\ Chen}, {A nonparametric view of network models and Newman--Girvan and other modularities},
  \emph{Proc. Nat.\ Acad.\ Sci.}
  \textbf{106}\,(50) (2009), {21068--21073}.
 


\bibitem{bickel2015correction}
        {P.~J.~Bickel,  A.~Chen,  Y.~Zhao, E.~Levina and J.~Zhu}, 
        {Correction to the proof of consistency of community detection},
        \emph{Ann.\ Statist.} \textbf{43} (2015), 462--466.


\bibitem{louvain}
	{V. D.~Blondel,  J.~Guillaume, R.~Lambiotte and E.~Lefebvre},
	{Fast unfolding of communities in large networks},
        \emph{J.\ Stat.\ Mech.\ Theory Exp.} 
        \textbf{2008}\,(10) (2008), P10008.  


\bibitem{bolla2015spectral}
	{M.~Bolla, B.~Bullins, S.~Chaturapruek, S.~Chen and K.~Friedl}, {Spectral properties of modularity matrices},	\emph{Linear Algebra Appl.} 
	\textbf{473}
	(2015), {359--376}.

\bibitem{BCLSVjournal}
	C. Borgs, J.~T.~Chayes, L.~Lov{\'a}sz, V.~T.~S{\'o}s and K.~Vesztergombi,
	{Convergent sequences of dense graphs {I}: {S}ubgraph frequencies, metric properties and testing},
	\textbf{219}\,(6) (2008), {1801--1851}. 

\bibitem{nphard}
	U.~Brandes, D.~Delling, M.~Gaertler,  R.~Gorke, M.~Hoefer, Z.~Nikoloski and D.~Wagner,  
        {On modularity clustering},
	\emph{Knowledge and Data Engineering, IEEE} 
	  \textbf{20}\,(2)
	(2008), 172--188.

\bibitem{chellig2022modularity}
        {J.~Chellig, N.~Fountoulakis and F.~Skerman},
        {The modularity of random graphs on the hyperbolic plane},
        \emph{J.\ Complex Networks}
        \textbf{10}\,(1)
        (2022). 




\bibitem{davis2018consistency}
        {E.~Davis and S.~Sethuraman},
        {Consistency of modularity clustering on random geometric graphs},
        \emph{Ann.\ Appl.\ Probab.}
        \textbf{28}\,(4)
        (2018), {2003--2062}.



\bibitem{demboMontanariSen}
	{A.~Dembo, A.~Montanari and S.~Sen},
    {Extremal cuts of sparse random graphs},
	\emph{Ann.\ Probab.}	
	\textbf{45}\,(2)
	(2017), {1190--1217}.


\bibitem{cookbook}
        {J.~D{\'\i}az and D.~Mitsche},
        {The cook-book approach to the differential equation method},
        \emph{Comp.\  Sci.\ Rev.}
        \textbf{4}\,(3)  (2010), {129--151}.


    
\bibitem{dinh2011finding}
	{T.~N.~Dinh and M.~T.~Thai},
	{Finding community structure with performance guarantees in scale-free networks},
	\emph{Privacy, Security, Risk and Trust (PASSAT) and 2011 IEEE Third Internat.\ Conf.\ on Social Computing (SocialCom)},
	{IEEE}
	(2011), {888--891}.


\bibitem{elgesem}
D.~Elgesem, L.~Steskal and N.~Diakopoulos, {Structure and
content of the discourse on climate change in the blogosphere: The big picture}, \emph{Environmental
Communication} (2014), {21--40}.


\bibitem{ferber2022friendly}
  {A.~Ferber, M.~Kwan, B.~Narayanan, A.~Sah  and M.~Sawhney},
  {Friendly bisections of random graphs},
  \emph{Comm.\ Amer.\ Math.\  Soc.}
  \textbf{210}  
  (2022), {380--416}.  

\bibitem{FortBart2008}
	S.~Fortunato and M.~Barth\'elemy,
        {Resolution limit in community detection},
	\emph{Proc.\ Nat.\ Acad.\ Sci.}
	\textbf{104}\,(1)
	(2007), {36--41}.

\bibitem{fortunato2016community}
	{S.~Fortunato and D.~Hric}, 
        {Community detection in networks: A user guide},
	\emph{Phys.\  Rep.}
	\textbf{659}
	(2016), {1--44}.
	
\bibitem{fountoulakis2021clustering}
        {N.~Fountoulakis, P.~Van Der Hoorn,  T.~M{\"u}ller and M.~Schepers},
        {Clustering in a hyperbolic model of complex networks},
        \emph{Electron.\ J.\ Probab.}
        \textbf{26}
        (2021), {1--132}.

\bibitem{friedman2008proof}
	J.~Friedman, 
        {A proof of Alon's second eigenvalue conjecture and related problems},
	\emph{Mem.\ Amer.\ Math.\ Soc.}            
	\textbf{195}\,(910)
	(2008), viii--100.

    


\bibitem{purplebook}
        {S.~Janson, T.~Luczak and A.~Rucinski}, {\em Random Graphs},
        Wiley, 2011.






\bibitem{kloks1994treewidth}
        T.~Kloks,
        {\em Treewidth: Computations and Approximations},
        Springer, 1994.
        
    


\bibitem{koshelev2023modularity}
        M.~Koshelev,
        Modularity in planted partition model,
        \emph{Computational Management Science}
        \textbf{20}\,(34)
        (2023). 


\bibitem{krioukov2010hyperbolic}
        D.~Krioukov, F.~Papadopoulos, M.~Kitsak, A.~Vahdat and M.~Bogun{\'a},
        {Hyperbolic geometry of complex networks},
        \emph{Phys.\ Rev.\ E}
        \textbf{82}\,(3) {036106}
        (2010).


\bibitem{lambiotte2021modularity}
        {R.~Lambiotte and M.~T.~Schaub},
        \emph{Modularity and Dynamics on Complex Networks},
        Cambridge Univ.\ Press,
        2021.


\bibitem{popular}
	{A.~Lancichinetti and S.~Fortunato},
	{Limits of modularity maximization in community detection},
	\emph{Phys.\ Rev.\ E.}
	\textbf{84}\,(6) {066122}
	(2011).

    
\bibitem{lasonsulkowska2023modularity}
        {M.~Laso{\'n} and M.~Sulkowska},
        {Modularity of minor-free graphs},
        \emph{J.\ Graph Theory} 
        \textbf{102}\,(4)   
        (2023), {728--736}.
  

\bibitem{lichev2022modularity}
        {L.~Lichev and D.~Mitsche},
        {On the modularity of 3-regular random graphs and random graphs with given degree sequences},
        \emph{Random Structures Algorithms}
        \textbf{61}\,(4)
        (2022), {754--802}.


\bibitem{liebenau2024asymptotic}
        {A.~Liebenau and N.~Wormald},
        {Asymptotic enumeration of graphs by degree sequence, and the degree sequence of a random graph},
        \emph{J.\ Europ.\ Math.\ Soc.}
        \textbf{26}\,(1)
        (2024), {1--40}.


\bibitem{modexpansion}
        {B.~Louf, C.~Mc{D}iarmid and F.~Skerman},
        {Modularity and graph expansion},
        {15th Innovations in Theoretical Computer Science Conference (ITCS 2024)}
        (2024), {78:1--78:21}.

\bibitem{lovasz2012large}
        {L.~Lov{\'a}sz},
        {\em Large Networks and Graph Limits},
        Vol.\ 60,
        \emph{Amer.\ Math.\ Soc.\ Colloq.\ Publ.},
        2012.

\bibitem{lusseau2003emergent} 
    D.~Lusseau,	
    {The emergent properties of a dolphin social network},
	\emph{Proc.\ Roy.\ Soc.\ London, Series B: Biological Sciences}
	\textbf{270}
	(2003), {S186--S188}.




\bibitem{mcdiarmid2023random}
        {C.~Mc{D}iarmid and S.~Saller},
        {Random graphs embeddable in order-dependent surfaces},
        \emph{Random Structures Algorithms}
        \textbf{64}\,(4)
        (2023), {1--46}.


\bibitem{treelike}
	C.~Mc{D}iarmid and F.~Skerman,
	{Modularity of regular and treelike graphs},
        \emph{J.\ Complex Networks}
	\textbf{6}\,(4)
	(2018), {596--619}.

\bibitem{ERmod}
	C.~Mc{D}iarmid and F.~Skerman,
        {Modularity of {E}rd{\H{o}}s--{R}{\'e}nyi random graphs,}
	\emph{Random Structures Algorithms}
	\textbf{57}\,(1)
	(2020), 211--243.


\bibitem{sampling}
	C.~Mc{D}iarmid and F.~Skerman,
	{Modularity and partially observed graphs},
	{arXiv:2112.13190},
	(2021).

\bibitem{vdense}
	C.~Mc{D}iarmid and F.~Skerman,
        {Modularity of nearly complete graphs and bipartite graphs},
        {arXiv:2311.06875},
	(2023).
	

\bibitem{majstorovic2014note}
	{S.~Majstorovic and D.~Stevanovic},
	{A note on graphs whose largest eigenvalues of the modularity matrix equals zero},
	\emph{Electron.\ J.\ Linear Algebra}
	\textbf{27}\,(1)
	(2014), {611--618}. 

\bibitem{w1hard}
	{K.~Meeks and F.~Skerman},
	{The Parameterised Complexity of Computing the Maximum Modularity of a Graph},
        \emph{Algorithmica},
	\textbf{82}\,(8)
	(2020), {2174--2199}.



  

\bibitem{distinguish}
	{E.~Mossel,  J.~Neeman and  A.~Sly},
	{Reconstruction and estimation in the planted partition model},
        \emph{Probab.\ Theory \& Related Fields}
	\textbf{162} 
        (2015), {431--461}.

\bibitem{NewmanGirvan}
	{M.~E.~J.~{N}ewman and M.~Girvan},
	{Finding and evaluating community structure in networks},
	\emph{Phys.\ Rev.\ E.}
	\textbf{69}\,(2)
	(2004), {026113}.
	
    


\bibitem{porter2009communities}
	{M.~A.~Porter, J.~Onnela and P.~J.~Mucha},
	{Communities in networks},
	\emph{Notices\ Amer.\ Math.\ Soc.}
	\textbf{56}\,(9)
	(2009), {1082--1097}.

\bibitem{prokhorenkova2017modularity}
        {L.~O.~Prokhorenkova, A.~Raigorodskii and P.~Pra{\l}at},
        {Modularity of complex networks models},
        \emph{Internet Math.}
        (2017), {115--126}.

\bibitem{trulymodular}
	{J.~Reichardt and S.~Bornholdt},
	{When are networks truly modular?},
	\emph{Phys.\ D.} 
	\textbf{224}\,(1)
	(2006), {20--26}.

\bibitem{robinson1992almost}
        {R.~Robinson and N.~C.~Wormald},
        {Almost all cubic graphs are Hamiltonian},
        \emph{Random Structures Algorithms}
        \textbf{3}\,(2)
        (1992), {117--125}.

\bibitem{rybarczyksulkowska2025PA}
        {K.~Rybarczyk and M.~Sulkowska},
        {Modularity of preferential attachment graphs},
        {arXiv:2501.06771}, (2025).

\bibitem{rybarczyksulkowska2025Gnp}
        {K.~Rybarczyk and M.~Sulkowska},
        {New bounds on the modularity of $G(n,p)$},
        {arXiv:2504.16254}, (2025).

\bibitem{thesis}
        {F.~Skerman},
        \emph{Modularity of Networks},
	{D.Phil.~Thesis},
        {University of Oxford},
	2016.


\bibitem{traag2019leiden}
        {V.A.~Traag, L.~Waltman and N.~J.~Van Eck},
        {From {L}ouvain to {L}eiden: guaranteeing well-connected communities},
        \emph{Scientific reports}
        \textbf{9}\,(1)
        (2019), {1--12}.

\bibitem{zhao2012consistency}  
        {Y.~Zhao, E.~Levina and J.~Zhu},
        {Consistency of community detection in networks under degree-corrected stochastic block models},
        \emph{Ann.\ Statist.}
        \textbf{40}\,(4)
        (2012), {2266--2292}.



\end{thebibliography}

\end{document}